\documentclass[runningheads]{llncs}
\usepackage{graphicx}
\usepackage{amsmath,amsfonts,amssymb,bbm}
\newcommand{\eop}{\hfill{$\blacksquare$}}

\newcommand{\up}{\Upsilon}
\usepackage[dvipsnames]{xcolor}
\usepackage{nicefrac}
\usepackage{enumerate,lipsum}
\usepackage[shortlabels]{enumitem}
\usepackage{multirow}
\usepackage{float}
\newcommand{\indc}[1]{\mathbbm{1}_{\{#1\}}}

\newcommand{\htheta}{\hat{\theta}}
\newcommand{\hpsi}{\hat{\psi}}

\newcommand{\hbeta}{\hat{\beta}}
\newcommand{\heta}{\hat{\eta}}
\newcommand{\hpi}{\pi(\hat{\beta})}

\usepackage{hyperref}

\newcommand{\hide}[1]{}

\newcommand{\TR}[2]{#2}  

\definecolor{blizblue}{rgb}{0.67, 0.9, 0.93}
\begin{document}
\title{Evolutionary Vaccination  Games with  premature vaccines   to combat  ongoing deadly pandemic } 

%
%
\author{Vartika Singh
\and
Khushboo Agarwal 
\and
Shubham
\and
Veeraruna Kavitha\thanks{The work of first and second author is partially supported by Prime Minister's Research Fellowship (PMRF), India.}
}
\authorrunning{ V. Singh et al.}
%
\institute{
IEOR, IIT Bombay
\\
\email{\{vsvartika, agarwal.khushboo, 19I190009, vkavitha\}@iitb.ac.in}
}
\maketitle              
\vspace{-6mm}
\begin{abstract}
 We consider a vaccination game that results with the introduction of premature and possibly scarce vaccines  introduced in a desperate bid to combat the otherwise ravaging deadly pandemic. The response of unsure agents amid many uncertainties makes this game completely different from the previous studies.  We construct a framework that combines SIS epidemic model with a variety of dynamic  behavioral vaccination responses and demographic aspects. The response of each agent is  influenced by the vaccination hesitancy and urgency, which arise due to their personal belief about efficacy and side-effects of the vaccine, disease characteristics,  and relevant reported information (e.g.,  side-effects, disease statistics etc.). Based on such aspects, we identify  the responses that are stable against static mutations.  By analysing the attractors of the resulting ODEs, we observe interesting  patterns in the limiting state of the system under evolutionary stable (ES) strategies, as a function of various defining parameters. There are responses for which the disease is eradicated completely (at limiting state), but none  are stable against mutations. Also, vaccination abundance  results in higher infected fractions at ES limiting state, irrespective of the disease death rate.
 
\vspace{-2mm}
\end{abstract}
%
%
%

\vspace{-10mm} 
\section{Introduction}
The impact of pandemic in today's world is unquestionable and so is the need to analyse various related aspects.
There have been many disease outbreaks in the past
 and the most recent Covid-19 pandemic is still going on. 
Vaccination is known to be of great help
; however, the effectiveness depends on the responses of the population (\cite{efv1} and references therein).
The vaccination process gets more challenging when vaccines have to be introduced prematurely, without much information about their efficacy or side effects, to combat the ravaging on-going  pandemic.  The scarcity of vaccines makes it all the more challenging. 


 There is a vast literature developed during the current pandemic that  majorly focuses on exhaustive experiments.  Recently authors in \cite{vg5}  discuss  the importance of game theoretic and social network models  for better understanding  the pandemic.  Our paper exactly  aims to achieve this purpose. We aim to develop a mathematical framework that mimics the  ongoing pandemic as closely as possible. \textit{We consider vaccination insufficiency,  hesitancy, impact of individual vaccination responses, possibility of 
 excess deaths, lack of information (e.g., possible end of the disease, vaccine details) etc. Our model brings together the well known  epidemic SIS  model (\cite{SIS1,SIS2}), evolutionary game theoretic framework (\cite{webb}), dynamic behavioural patterns of the individuals and demographic aspects. }

Majority of the literature on vaccination games assumes some knowledge which aids in vaccination plans;  either they consider seasonal variations, or the time duration and or the time of occurrence of the disease is known  etc.  For example, in \cite{vg1,vg3} authors consider a replicator dynamics-based vaccination game; each round occurs in two phases, vaccination and the disease phase. Our paper deals with agents, that continually operate under two contrasting fears,  vaccination fear
and the deadly pandemic fear. 
They choose between the two fears, by estimating the perceived   cost of the two aspects.  
 In \cite{vg2} authors consider perceived vaccination costs as in our model and study the `wait and see' equilibrium. Here the agents  choose one among 53 weeks to get vaccinated. In contrast, {\it we do not have any such information}, and the aim is to eradicate the  deadly disease.
 
 \hide{
 {\color{red}In \cite{vg4}, authors focus on the impact of lack of information and awareness on vaccination program. 
The authors in \cite{vg6} study childhood vaccination program as decided by the parents. 
 These papers focus on very different aspects related to vaccination games. 

Follow-the-crowd behavior is seen in replicator dynamics, where individuals learn from the experiences of their neighbors. However, in our model, follow-the-crowd agents learn from the reported information.}}

In particular, we consider a population facing an epidemic with limited availability of vaccines. 
Any individual can  either be infected, susceptible or vaccinated.\textit{ It is evident nowadays (with respect to Covid-19) that recovery does not always result in immunization. 
 By drawing parallels, the recovered individuals become  susceptible again.}   We believe the individuals in this current pandemic are divided in opinion due to vaccination hesitancy and vaccination urgency. A vaccination urgency can  be observed depending upon the reported information about the sufferings and deaths due to disease, lack of hospital services etc. On the other hand, the vaccination hesitancy could be due to individuals' belief in the efficacy  and the emerging fears due to the  reported side effects. 

To capture above factors, we model a  variety of possible vaccination responses. We first consider the individuals who exhibit the \textit{follow-the-crowd (FC)} behavior, i.e., their confidence (and hence inclination) for vaccination  increases as more individuals get vaccinated. 
In other variant, the interest of an individual in vaccination reduces as vaccinated proportion further grows. Basically such individuals attempt to enjoy the resultant benefits of not choosing vaccination and still being prevented from the disease; we refer them as \textit{free-riding (FR) agents}. Lastly, we consider individuals who make more informed decisions based also on the infected proportion; these agents exhibit  increased urge towards vaccination as the infected population grows. We name such agents as \textit{vigilant agents}.

In the era of pre-mature vaccine introduction and minimal information, the agents depend heavily on the perceived cost of the two contrasting (vaccination/infection cost) factors.  The  perceived cost of infection can be large leading to a vaccination urgency, when there is vaccine scarcity. In case there is abundance, one may perceive a smaller risk of infection  and procrastinate the vaccination till the next available opportunity. \textit{Our results interestingly indicate a larger infected proportion at ES equilibrium as vaccine availability improves.}



We intend to investigate the resultant of the vaccination response of the population and the nature of the disease. Basically, we want to understand if the disease can be overpowered by given type of vaccination participants. Our analysis is layered: firstly, {\it we use stochastic approximation techniques to derive  time-asymptotic proportions and identify all possible equilibrium states, for any given vaccination response and disease characteristics}.  Secondly, we derive the vaccination responses that are \textit{stable against mutations}. We  slightly modify the definition of classical evolutionary stable strategy, which we refer as \textit{evolutionary stable strategy against static mutations (ESS-AS)}\footnote{Agents use dynamic policies, while mutants use static variants (refer Section \ref{sec_prob_desc}).}. We study the equilibrium/limiting states reached by the system under ES vaccination responses.  

Under various ES (evolutionary stable) vaccination responses, the dynamic behaviour at the beginning could be different, but \textit{after reaching equilibrium, individuals either vaccinate with probability one or zero}. Some interesting patterns are observed in ES limiting proportions as a function of important defining parameters. For example, the ES limit infected proportions are concave functions of birth rate
. \textit{With increased excess deaths, we have smaller infected proportions at ES limiting state}. Further, \textit{there are many vaccination responses which eradicate the disease completely at equilibrium; but none of them are evolutionary stable, unless the disease can be eradicated without vaccination}.
At last, we corroborate our theoretical results   by performing exhaustive Monte-Carlo simulations. 

\vspace{-6mm}
\section{Problem description and background}\label{sec_prob_desc}

\vspace{-2mm}
We consider a  population facing an epidemic, where at time $t$, the state of the system  is {\small$(N(t),S(t),V(t),I(t))$}. These    respectively represent   total, susceptible, vaccinated and  infected population. Observe {\small$N(t) =S(t)+ V(t)+I(t) $}.

At any time $t$, any susceptible individual can contact anyone among the infected population according to exponential distribution with parameter $\lambda/N(t)$, $\lambda > 0$ (as is usually done in epidemic models, e.g., \cite{SIR_proof}). In particular, a contact between a susceptible and an infected individual may result in spread of the disease, which is captured by  $\lambda$. Any infected individual may recover after an exponentially distributed time with parameter $r$ to  become susceptible again. A susceptible individual can think of vaccination  after exponentially distributed time with parameter $\nu$. At that epoch, the final decision (to get vaccinated) depends on the  information available to the individual. We refer to the probability of a typical individual getting vaccinated as $q$; more details will follow below.
{\it It is important to observe here that 
$\nu$ will also be governed by vaccination availability;  the individual decision rate is upper bounded by availability rate. }
Further, there can be a birth after exponentially distributed time with parameter $b\hspace{0.5mm} N(t)$. A death is possible in any compartment after exponentially distributed time with parameter $ d   \hspace{0.7mm} I(t)$ (or $d\hspace{0.7mm} S(t)$ or $d \hspace{0.7mm} V(t)$), and excess death among infected population  with parameter $d_e \hspace{0.7mm} I(t)$.  Furthermore, we assume $b>d+d_e$.

One of the objectives is to analyse the (time) asymptotic proportions of the infected, vaccinated and susceptible population depending upon the disease characteristics and vaccination responses. To this end, we consider the following fractions: \vspace{-8mm}

{\small 
\begin{align}\label{eqn_theta_t}
\hspace{1mm}\theta(t):=\frac{I(t)}{N(t)},\ \psi(t):=\frac{V(t)}{N(t)}, \text{ and } \phi(t):=\frac{S(t)}{N(t)}.
\end{align}}

\vspace{-4mm}
\subsection{Responses towards vaccination }\label{subsec_agents_bhvr}

In reality, the response towards available vaccines depends upon the cost of vaccination, the severity of infection and  related recovery and death rates. However, with lack of information, the vaccination decision of any susceptible   depends heavily on the available data,  ($\theta(t)$, $\psi(t)$), the  fractions of infected and vaccinated population at that time. Further, it could depend on the willingness of the individuals towards vaccination. The final probability of getting vaccinated at the decision/availability epoch is captured by  probability  $q=q(\theta,\psi,\beta)$; this probability is also influenced by the parameter $\beta$ which is a characteristic of the given population. Now, we proceed to describe different possible behaviours  exhibited by agents/individuals in the population:

{\bf Follow the crowd (FC) agents:}
These type of agents make their decision to get vaccinated by considering only the vaccinated proportion of the population; they usually ignore other statistics. They believe vaccination in general is good, but are hesitant because of the possible side effects. The hesitation reduces as the proportion of vaccinated population increases, and thus  accentuating the likelihood of vaccination for any individual. In this case, the probability that any individual will decide to vaccinate is given by  $q := \beta \psi(t)$. 

{\bf Free-riding (FR) agents:}
\hide{
The agents can also exhibit a more evolved behaviour than Follow the crowd, in particular free-riding behaviour. Here, the agents decide to get vaccinated on the basis of the fraction of vaccinated population, just like FC agents. But at the same time, they prefer to benefit from the collective decisions of the others (to get vaccinated) than going for vaccination themselves, thereby preventing themselves from the risks associated with vaccine or paying for its cost. We call such agents are free-riding agents.

These agents in order to avoid any risks related to vaccinations, observe the crowd behaviour and act the opposite. Their motive is to survive without accumulating the risks related to vaccination and  behave as a free rider.

These agents tend to go for vaccination as the proportion of vaccinated population increases, but they don't just follow the crowd. Instead, as the proportion of vaccinated population grows, their tendency to vaccinate starts to reduce. They rely free benefits of vaccination of others, and hence, we name them `Free riding agents'.
}
These agents exhibit a more evolved behaviour than FC individuals, in particular free-riding behaviour. In order to avoid any risks related to vaccination or paying for its cost, they observe the crowd behaviour more closely. Their willingness towards vaccination also improves with $\psi(t)$, however they tend to become free-riders beyond a limit. In other words, as $\psi(t)$ further grows, their tendency to vaccinate starts to diminish. We capture such a behaviour by modelling the probability  as $q:= \beta \psi(t)(1-\psi(t)).$

{\bf Vigilant agents:}
Many a times individuals are more vigilant, and might consider vaccination only if the infected proportion  is above a certain threshold, in which case we set  $q= \beta  \psi(t) \indc{\theta (t) > \Gamma}$ ($\Gamma $ is an appropriate constant). This dependency could also be continuous, we then model the probability as $q= \beta \theta(t) \psi(t)$. We refer to them as Vigilant Follow-the-Crowd (VFC) agents. 

\vspace{2mm}
\noindent{\bf Vaccination Policies:} 
We refer the vaccination responses of the agents as vaccination policy, which we represent by $\pi(\beta)$. When the vaccination response policy is $\pi(\beta)$,   the agents get  vaccinated with  probability $q_{\pi(\beta)}(\theta, \psi)$ at the vaccination decision epoch,  based on the system state $(\theta, \psi)$ at that epoch. For ease of notation, we avoid $\pi$ and $\beta$ while representing $q$. 

As already explained,  we consider the following vaccination policies: a) when $\pi(\beta) = F_C (\beta)$, i.e, with follow-the-crowd policy, the agents choose  to vaccinate with probability \underline{ $q(\theta, \psi) = \min\{1, \tilde{q} (\theta, \psi) \}, \  \tilde{q} (\theta, \psi) =  {\beta} \psi$} at vaccination decision epoch;  b) when $\pi (\beta) = F_R(\beta)$, the free-riding policy,  $\tilde{q}(\theta, \psi) =  
{\beta} \psi (1-\psi)$; and (c)
 the policy $\pi (\beta)  = V_{FC}^1 (\beta)$ represents vigilant (w.r.t. $\theta$) follow the crowd (w.r.t. $\psi$) policy 
and $\tilde{q}(\theta, \psi) = {\beta} \theta \psi$. We define $\Pi:=\{ F_C,F_R, V_{FC}^1\}$ to be the set of all these policies. We discuss the fourth type of policy  $V_{FC}^2$ with $\tilde{q}(\theta,\psi)=\beta\theta\indc{\theta>\Gamma}$ separately, as the system responds very differently to these agents. Other   behaviour patterns are for future study.

\vspace{-2mm}
\subsection{Evolutionary behaviour}

The aim of this study is three fold: (i) to study the dynamics and understand the equilibrium states (settling points) of the disease depending upon the agents' behavior and the availability of the vaccine, (ii) to compare and understand the differences in the  equilibrium states depending on agents' response towards vaccine, and (iii) to investigate if these equilibrium states are stable against mutations, using the well-known concepts of evolutionary game theory. Say a mutant population of small size invades such a system in equilibrium. We are interested to investigate if the agents using the original (vaccination response) are still better than the mutants. Basically we use the concept of evolutionary stable strategy (ESS), which in generality is defined as below (e.g., \cite{webb}):

By definition, for $\pi({\beta})$ to be an ESS, it should satisfy two conditions, i) $\pi({\beta})\in \arg\min_{\pi \in \widetilde{\Pi}} u(\pi ; \pi({\beta}))$, where $\widetilde{\Pi}$ is the set of the policies, $u(\pi, \pi')$ is the utility/cost of the (can be mutant) user that adapts policy $\pi$, while the rest of the population uses policy $\pi'$;  and ii) it should be stable against mutations, i.e. there exist an $\bar{\epsilon}$, such that for all {$\epsilon\le\bar{\epsilon}  ( \pi)$},
$$u(\pi,\pi_\epsilon({\beta},\pi))>u(\pi({\beta}),\pi_\epsilon({\beta},\pi)) \mbox{  for any } \pi \ne \pi(\beta),$$
where $\pi_\epsilon({\beta},\pi):= \epsilon \pi +(1-\epsilon)\pi({\beta})$, represents the policy in which $\epsilon$ fraction of agents (mutants) use strategy $\pi$ and the other fraction uses $\pi({\beta})$. 

In the current paper, we restrict our definition of ESS to cater for mutants that use static policies, $\Pi^D$.  Under any static policy $\pi\equiv q$, the agent gets vaccinated with constant probability $q$ at any  decision epoch, irrespective of the system state.  
We now \textit{define  the Evolutionary Stable Strategies, stable  Against Static mutations} and the exact definition is given below.  

\begin{definition}{\bf[ESS-AS]} A policy $\pi(\beta)$ is said to be ESS-AS,  
i) if $\{q^*_{\pi(\beta)} \} = {\cal B} (\pi(\beta) )$, where the static-best response set  \vspace{-2mm}
$${\cal B} (\pi(\beta) ) := \arg \min_{q \in [0,1]}  u(q, \pi(\beta) ),
$$ and  $q^*_{\pi(\beta)} = q (\theta^*,\psi^*)$  
is the probability with which the agents get vaccinated after the system reaches equilibrium under strategy $\pi(\beta)$;  and 
ii) there exists an $\bar{\epsilon}$ such that $\{q^*_{\pi(\beta)} \} = {\cal B} (\pi_\epsilon (\beta, q) )  $, for any $\epsilon \le {\bar \epsilon}(q)$ and any $q$. 
\end{definition}

\noindent{\bf Anticipated utility/cost of a user:}
Once  the population settles to an equilibrium (call it $\hat{\theta}
, \hat{\psi}$) under  a certain policy $\pi(\hat{\beta})$,  the users (assume to) incur a  cost  that  depends upon various factors. To be more specific, any user estimates\footnote{We assume mutants are more rational, estimate various rates using reported data.
} its overall 
cost of vaccination considering the pros and cons as below to make a judgement about vaccination. 

The cost of infection (as perceived by the user)  can be summarized by   $p_I(\hat{\theta})(c_{I_1} +  c_{I_2}d_e \hat{\theta})$, where  $p_I(\hat{\theta})$ equals the probability that the user gets infected before its next decision epoch (which depends upon the fraction of infected population $\hat{\theta}$ and availability/decision rate $\nu$),  $c_{I_1}$ is the cost of infection without death (accounts for the sufferings due to disease, can depend on $r$),  while $c_{I_2}d_e \hat{\theta}$ is the perceived chance of death after infection. Observe here that $d_e \hat{\theta}$ is the fraction of excess deaths among infected population which aids in this perception. 

On the contrary, 
the cost of vaccination is summarized by  $c_{v_1} + \min\left\{\bar{c}_{v_2}, \nicefrac{c_{v_2}}{\hat{\psi}}\right\}$, where  $c_{v_1}$ is the  actual cost of vaccine. 
Depending upon the fraction $\hat{\psi}$ of the population vaccinated and their experiences, the users   anticipate additional cost of vaccination (caused due to side-effects) as captured by the second term $\nicefrac{c_{v_2}}{\hat{\psi}}$. Inherently {\it we assume here that the side effects are not significant,} and hence in a system  with a bigger vaccinated fraction, the vaccination hesitancy is lesser. Here $\bar{c}_{v_2}$  accounts  for maximum  hesitancy.   In all, the expected anticipated cost of vaccination by a user  in a  system at equilibrium (reached under   $\pi(\hat{\beta)}$)  equals: 
the probability of vaccination  (say $q$) times the anticipated cost of vaccination, plus  $(1-q)$ times the anticipated cost of infection. Thus we define:

\begin{definition}{\bf[User utility at equilibrium]}

When the population is using policy $\pi(\hat{\beta})$ and an agent attempts to get  itself vaccinated with probability $q$, then,  the user utility function is given by:
\begin{eqnarray}\label{eqn_utility}
  u(q; \pi(\hat{\beta})) &:=& q\left(c_{v_1}+\min\left\{\bar{c}_{v_2}, \frac{c_{v_2}}{\hat{\psi}}\right\}\right)+(1-q)p_I(\hat{\theta})(c_{I_1} +  c_{I_2}d_e \hat{\theta})\nonumber\\
  &=& q h(\pi(\hat{\beta}))  + p_I(\hat{\theta})(c_{I_1} +  c_{I_2}d_e \hat{\theta}), \mbox{ where }\\
  h(\pi(\hat{\beta}))&=& h(\hat{\theta}, \hat{\psi}) := c_{v_1} + \min\left\{\bar{c}_{v_2}, \frac{c_{v_2}}{\hat{\psi}} \right\} - p_I(\hat{\theta})(c_{I_1} +  c_{I_2}d_e \hat{\theta}).\nonumber
\end{eqnarray}

\end{definition}


In the next section, we begin with ODE approximation of the system, which facilitates in deriving the limiting behaviour of the system. Once the limiting behaviour is understood, we proceed  towards evolutionary stable strategies.

\section{Dynamics and ODE approximation}\label{ODE_approx}
 Our aim in this section is to understand the limiting behaviour of the given system.
The system is transient with $b > d+d_e$; it is evident that the population would not settle to a stable distribution (it would explode as time progresses with high probability). However the fraction of people in various compartments (given by \eqref{eqn_theta_t})  can possibly reach some equilibrium and we look out for this equilibrium or limiting proportions (as is usually considered in literature (\cite{ltp})). 
 
 To study the limit proportions, it is sufficient to analyse the process at transition epochs. Let $\tau_k$ be the $k^{th}$ transition epoch, and infected population immediately after $\tau_k$ equals   $I_k := I(\tau_k^+) = \lim_{t \downarrow \tau_k} I(t)$; similarly, define $N_k, S_k$ and $V_k$. Observe here that $T_{k+1} := \tau_{k+1} - \tau_k$ is exponentially distributed with a parameter that depends upon previous system state $(N_{k}, S_{k}, I_{k}, V_k)$.

\textbf{Transitions:} Our aim is to derive the (time) asymptotic fractions of \eqref{eqn_theta_t}. Towards this, we define the same fractions at transition epochs, 
$$\theta_k:=\frac{I_k}{N_k}, \psi_k:=\frac{V_k}{N_k}, \text{ and } \phi_k:=\frac{S_k}{N_k}. $$ Observe that   $\theta_k+\psi_k+\phi_k=1$. To facilitate our analysis, we also define a slightly different  fraction,  $\eta_k:=N_k/k$ for $k > 1$ and  $\eta_0 := N(0)$, $\eta_1 := N(1)$.
As described in previous section, the size of the infected population evolves between two transition epochs according to:
\begin{align}
\label{Eqn_Ik}
        I_{k+1} & = I_k + G_{I, k+1}, \mbox{ with }  G_{I, k+1} := {\mathbb I}_{k+1} - {\mathbb R}_{k+1}- {\mathbb D}_{I,k+1},
        \end{align}
        where ${\mathbb I}_{k+1}$ is the indicator that the current epoch is due to a new infection,  ${\mathbb R}_{k+1}$ is  the indicator of a recovery and  ${\mathbb D}_{I,k+1}$ is the indicator that the current epoch is due to a death among the infected population. 
        Let ${\cal F}_k := \sigma  (I_j, S_j, V_j, N_j, j \le k)$ represent the sigma algebra generated by the history until the observation epoch $k$ and let $E_k [\cdot]$  represent the corresponding  conditional expectation.
       By conditioning on ${\cal F}_k$, using the memory-less property of exponential random variables,
       \begin{eqnarray}\label{eqn_cond_exp}
            E_k [{\mathbb I}_{k+1}] 
            &=& \frac{\frac{\lambda I_k S_k}{N_k}}{ N_k b + N_k d + I_k d_e + \frac{\lambda I_k S_k}{N_k}+ S_k\nu  +  I_k r }\ = \ \frac{\lambda \theta_k \phi_k }{ \varrho_k }, \mbox{ with, }\nonumber\\
        \varrho_k    & :=
             &   b  + d +d_e \theta_k +  \lambda \theta_k \phi_k+ \nu \phi_k + r \theta_k, \nonumber\\
               E_k [{\mathbb R}_{k+1}] 
            & =  &
             \frac{ r \theta_k }{\varrho_k } , \mbox{ and, }
            E_k [{\mathbb D}_{I,k+1}] 
            \ =\   
              \frac{ \theta_k  (d+d_e)}{\varrho_k }
          .
       \end{eqnarray}
         In similar lines the remaining types of the population evolve according to the following, where $\mathbb{V}_{k+1}$,  $\mathbb{B}_{k+1}$, $\mathbb
            {D}_{V,k+1}$ and ${\mathbb D}_{S,k+1}$ are respectively the indicators of vaccination, birth and corresponding deaths,  
\begin{align}
\label{Eqn_VNSk}
            V_{k+1} &= V_k \ + G_{V, k+1},  \ G_{V, k+1} := \mathbb{V}_{k+1}- \mathbb
            {D}_{V,k+1},\\
    N_{k+1}     &= N_k+ G_{N,k+1},  \ G_{N,k+1}:= \mathbb{B}_{k+1} -\mathbb
            {D}_{I,k+1}-\mathbb
            {D}_{V,k+1}-\mathbb
            {D}_{S,k+1},
    \mbox{ and,}
    \nonumber
    \\
    S_{k+1} &= N_{k+1} - I_{k+1} - V_{k+1}.    \nonumber
\end{align} 
As before, 

\vspace{-9mm}
 {\small      \begin{eqnarray}\label{eqn_cond_exp_2}
            E_k [{\mathbb V}_{k+1}] 
            &=&   \frac{\nu q(\theta_k, \psi_k) \phi_k }{\varrho_k },  \  \ 
            E_k [{\mathbb D}_{V,k+1}] 
            \ = \  
              \frac{ \psi_k d }{\varrho_k},\nonumber \\  
              E_k [{\mathbb D}_{S,k+1}]  &=&              \frac{ \phi_k d }{\varrho_k} , \mbox{ and, }
            E_k [{\mathbb B}_{k+1}] 
            =
             \frac{ b}{\varrho_k }.
       \end{eqnarray}}
  Let $\Upsilon_k:=[\theta_k,\psi_k,\eta_k]^T$. The evolution of $\Upsilon_k$ can be studied by a three dimensional system, described in following paragraphs.  To facilitate tractable mathematical analysis we consider a slightly modified system that freezes once $\eta_k $ reaches below a fixed small constant $\delta >0$. The rationale and the justification behind this modification is two fold: a) once the population reaches below a significantly small threshold, it is very unlikely that it explodes and the limit proportions in such paths are  no more interesting; b)  the initial population  $N(0)$ is usually large, let {\small $\delta ={2}/({N(0)-1})$} and then with $b > d+d_e$, it is easy to verify that  the probability, {\small $P\Big (\eta_k < \delta \mbox{ for some } k\Big |N(0) \Big ) \to 0$} as $N(0) \to \infty$. 
  From  \eqref{Eqn_Ik}, 
  
  \vspace{-4mm}
  {\small
\begin{eqnarray}\label{eqn_stc_theta}
\hspace{8mm}
 \frac{I_{k+1}}{N_{k+1}} &=& \frac{I_k}{N_k} +  \frac{I_{k+1}}{N_{k+1}}- \frac{I_k}{N_k}=\frac{I_k}{N_k} + \frac{1}{k+1}\frac{k+1}{N_{k+1}}\left[  I_{k+1}- \frac{I_k N_{k+1}}{N_k}\right],\nonumber\\
 &&
 \hspace{-10mm} = \  \frac{I_k}{N_k} + \frac{1}{k+1}\frac{k+1}{N_{k+1}}\left[ G_{I,k+1} -   \frac{N_{k+1}- N_k}{N_{k}} I_k  \right], \mbox{ \normalsize and thus including $\indc{\eta_k > \delta}$, }\nonumber\\
\theta_{k+1}    &=& \theta_k + \epsilon_k\frac{\indc{\eta_k > \delta}}{\eta_{k+1}}\left[ G_{I,k+1} -   (N_{k+1}-N_k)\theta_k\right], \ \ \epsilon_k:=\frac{1}{k+1}.
\end{eqnarray}}We included the indicator $\indc{\eta_k > \delta}$, as none of the population types change (nor there is any  evolution) once the population gets almost extinct. 
Similarly, from equation \eqref{Eqn_VNSk}, \vspace{-5mm}
\begin{eqnarray}\label{eqn_stoch_approx} \hspace{6mm}
\psi_{k+1} &=& \psi_k + \epsilon_k\frac{\indc{\eta_k >\delta}}{\eta_{k+1}}\left[ G_{V,k+1} - (N_{k+1}- N_k) \psi_k   \right],\\
    \eta_{k+1} &=& \eta_k +\indc{\eta_k > \delta} \epsilon_k\left[ G_{N,k+1}-\eta_k \right].\nonumber
\end{eqnarray} 
We  analyse this system using the results and techniques of  \cite{kushner2003stochastic}. In particular, we prove equicontinuity in extended sense for our non-smooth functions  (e.g., $q(\theta,\psi)$ may only be measurable), and then use \cite[Chapter 5, Theorem 2.2]{kushner2003stochastic}. Define $L_{k+1}:=[L^{\theta}_{k+1},L^{\psi}_{k+1},L^{\eta}_{k+1}]^T$, with \begin{eqnarray}\label{eqn_Lk}
     L^{\theta}_{k+1} &=&\frac{\indc{\eta_k > \delta}}{\eta_{k+1}}\left[ G_{I,k+1} -   (N_{k+1}- N_k) \theta_k  \right],\\
     L^{\psi}_{k+1}&=&\frac{\indc{\eta_k > \delta}}{\eta_{k+1}}\left[  G_{V,k+1}- (N_{k+1}- N_k) \psi_k   \right],\mbox{ and,  } L^{\eta}_{k+1}= \indc{\eta_k > \delta}(G_{N,k+1}-\eta_k).\nonumber
\end{eqnarray} 
Thus \eqref{eqn_stc_theta}-\eqref{eqn_stoch_approx}  can be rewritten as,      
 $\Upsilon_{k+1}=\Upsilon_k+\epsilon_k L_{k+1}$. 
Conditioning as in \eqref{eqn_cond_exp}:
\begin{eqnarray}\label{en_L_{k+1}}
        E_k[L^{\theta}_{k+1}]&=& \frac{ \theta_k\indc{\eta_k>\delta}}{\eta_k\varrho_k}\left[\phi_k  \lambda - r -d_e  -(b  - d_e \theta_k )\right]+ \alpha^{\theta}_k \\
    &=:&   g^\theta (\Upsilon_k)+ \alpha^{\theta}_k ,  \mbox{ where, }  \alpha^{\theta}_k = E_k \left [ L_{k+1}^\theta - \frac{\eta_{k+1} }{\eta_{k}} L_{k+1}^\theta  \right ]\nonumber.
\end{eqnarray}
In exactly similar lines,  we define 
 $g(\Upsilon_k)=[g^\theta (\Upsilon_k), g^\psi (\Upsilon_k), g^\eta (\Upsilon_k)]^T$  (details just below) and $\alpha^\psi_k$ such that $E_k[L^\psi_{k+1}]=g^\psi(\Upsilon_k)+ \alpha^\psi_k$ and  $E_k[L^\eta_{k+1}]=g^\eta(\Upsilon_k)$.
 Our claim is that the error terms
would converge to zero (shown by Lemma \ref{Lemma_alpha_km} in Appendix) and 
  ODE $\dot{\Upsilon} =g(\Upsilon)$ approximates the system dynamics,   where,  
\begin{align}\label{eqn_ODE}
\begin{aligned}
    g^\theta(\Upsilon) &= \frac{\theta\indc{\eta >\delta}}{\eta \varrho } \left[\phi\lambda - r - d_e- (b  - d_e \theta )\right], \  \phi = 1- \theta- \psi \\
  g^\psi(\Upsilon) &= \frac{\indc{\eta >\delta}}{\eta \varrho}\left[q(\theta,\psi)\phi \nu - (b - d_e \theta)\psi\right], \mbox{ and, }\\
   g^\eta (\Upsilon) &= \indc{\eta>\delta}\left(\frac{b-d - d_e \theta}{\varrho} - \eta\right), \  \ \varrho =  b  + d +d_e \theta +  \lambda \theta \phi+ \nu \phi + r \theta.
\end{aligned}
\end{align}

We now state our first main result (with proof in Appendix A).
\begin{enumerate}[{\bf A.}]
    \item Let the set $A$ be locally asymptotically stable in the sense of Liapunov for the ODE \eqref{eqn_ODE}.  Assume that $\{\up_n\}$ visits a compact set, $S_A$, in the domain of attraction, $D_A$, of $A$ infinitely often (i.o.) with probability  $\rho > 0$.
\end{enumerate}
\begin{theorem}\label{thrm1}
  Under assumption \textbf{A.}, i) the sequence converges, $\up_n \to A$ as $n \to \infty$ with probability at least $\rho$; and ii) for every $T>0$, almost surely there exists a sub-sequence $(k_m)$ such that: $(t_k:=\sum_{i=1}^{k}\epsilon_i)$
            $$
            \sup_{k: t_k \in [t_{k_m}, t_{k_m} + T]} d(\up_k, \up_{*}(t_k - t_{k_m})) \to 0,  \mbox{ as } m \to \infty, \mbox{ where, }
            $$
        $\up_*(\cdot)$ is the  solution of ODE \eqref{eqn_ODE} with initial condition  $\up_*(0) = \lim_{k_m} \up_{k_m}$. 
  \eop
\end{theorem}

\noindent Using above Theorem, one can derive the limiting state  of the system using that of the ODE (in non-extinction  sample paths, i.e., when $\eta_k > \delta$ for all $k$). 
Further, for any finite time window, there exists a sub-sequence along which  the disease dynamics are approximated by the solution of the ODE. The ODE should initiate at the limit of the system along such sub-sequence. 


\vspace{-4mm}
\section{Limit proportions and ODE attractors}

So far, we have proved that the embedded process of the system can be approximated by the solutions of the ODEs \eqref{eqn_ODE} (see Theorem \eqref{thrm1}). We will now analyse the ODEs and look for equilibrium states for a given  vaccination policy  $\pi(\hat{\beta})$. 
The following notations are used throughout:\textit{ we represent the parameter by $\hat{\beta}$ and the corresponding equilibrium states by $(\hat{\theta}, \hat{\psi})$. Let $\hat{q} : = q(\htheta, \hpsi)$ and $\hat{\varrho} := \varrho(\hat{\theta}, \hat{\psi})$.} In the next section, we identify the evolutionary stable (ES) equilibrium states $(\theta^*, \psi^*)$, among these equilibrium states, 
 and the corresponding vaccination policies $\pi(\beta^*)$.

We now identify the attractors  of the ODEs \eqref{eqn_ODE} that are locally asymptotically stable in the sense of Lyapunov (\textit{referred to as attractors}), which is a requirement of the assumption \textbf{A}. However, we are yet to identify the domains of attraction, which will be attempted in future. 
Not all infectious diseases lead to deaths. One can either have: (i) \textit{non-deadly disease}, where only natural deaths occur, $d_e = 0$, or (ii) \textit{deadly disease} where in addition, we have excess deaths due to disease,  $d_e > 0$.
We begin with the non-deadly case and  FC agents.  

The equilibrium states for  FC agents 
($\hbeta \geq 0$) are (proof in Appendix B):
\begin{theorem}{\bf [FC agents]} \label{thrm_FCagents}
Define $\rho := \lambda/(r+b+d_e)$, $\mu := b/\nu$.
When $d_e = 0$ and $\tilde{q}(\theta, \psi) = \hbeta \psi$,  at   the attractor we have $\hat{\eta} = \nicefrac{(b-d)}{{\hat \varrho}}$.
The remaining details of the attractors are in Table \ref{table_FC}.
\hide{
Further, 
(i) with $\hbeta < \rho \mu$ and $\rho > 1$ we have  $(\hat{\theta}, \hat{\psi}) =  \left( 1-\nicefrac{1}{\rho} , 0\right) $, (ii) they equal $(\hat{\theta}, \hat{\psi}) =  \left (0,   1- \nicefrac{\mu}{\hbeta} \right )$, if $\hbeta > \rho \mu$   and finally,  (iii)  we have, $(\hat{\theta}, \hat{\psi}) =  \left (0,0 \right )$ when $\rho < 1$ and $\mu > \hbeta$. }
The interior attractors $($when $(\htheta, \hpsi) \in (0,1)\times (0,1))$ are the zeros of the right hand side (RHS) of  ODE \eqref{eqn_ODE}.\eop
\end{theorem}

\renewcommand{\arraystretch}{1.5} 

\begin{table}[htbp]
\vspace{-4mm}
\centering
\begin{tabular}{|c|c|c|c|}
\hline
Nature  & \multicolumn{2}{c|}{Parameters}  & $(\htheta, \hpsi)$ \\ \hline \hline
\multirow{5}{*}{\begin{tabular}[c]{@{}c@{}} Endemic, \\ $\rho > 1$\end{tabular}}  & $\hbeta < \mu \rho$ &                                                        & $\left( 1-\frac{1}{\rho} , 0 \right) $ \\ \cline{2-4} 
      & \multirow{1}{*}{ $ \hbeta > \mu \rho$ , $\tilde{q}(0, 1-\frac{\mu}{\hbeta})  < 1 $ implies $\hbeta < \mu + 1$* }  & & $\left (0,   1- {\frac{ \mu}{\hbeta} }\right )$       \\ \cline{2-4}  
         &   \multirow{2}{*}{$ \hbeta > \mu \rho$ , $\tilde{q} (0, 1-\frac{\mu}{\hbeta})  > 1 $ implies $\hbeta > \mu + 1$}    & $\mu \rho < \mu + 1$ & $\left(0, \frac{1}{\mu+1} \right)$ \\ \cline{3-4}
 &   & $\mu \rho > \mu + 1$ & $\left(\theta_E, \psi_E \right )$ \\ \hline
SE, $\rho < 1$ &  $\mu > \hbeta$  &    & $(0, 0)$   \\ \hline
\end{tabular}
 \caption{Attractors for FC agents, $(\theta_E, \psi_E) := \left(1-\frac{1}{\rho}-\frac{1}{\mu \rho} ,
\ \frac{1} { \mu \rho} \right )$ \label{table_FC}}
\vspace{-4mm}
\end{table}

In all the tables of this section, the $*$ entries are also valid when $\rho < 1$. When $\hbeta = \mu \rho$,  
the ODE (and hence the system) is not stable; such notions are well understood in the literature and we avoid such marginal cases.  
%
%
We now consider the FR agents  (proof again in Appendix B):
\begin{theorem}{\bf [FR agents]}\label{thrm_FR_agents}
When   $\tilde{q}(\theta, \psi) = \hbeta \psi (1-\psi)$ and   $d_e = 0$, then the attractors for ODE \eqref{eqn_ODE} are $(\htheta, \hpsi, \nicefrac{(b-d)}{{\hat \varrho}})$, which are provided in Table \ref{table_FR}. The interior attractors are the zeros of the RHS of  ODE \eqref{eqn_ODE}. \eop
 \end{theorem}

\begin{table}[htbp]
\vspace{-4mm}
    \centering
\begin{tabular}{|c|c|c|c|}
\hline
Nature  & \multicolumn{2}{c|}{Parameters}  & $(\htheta, \hpsi)$ \\ \hline \hline
\multirow{5}{*}{\begin{tabular}[c]{@{}c@{}} Endemic, \\ $\rho > 1$\end{tabular}}  & $\hbeta < \mu \rho$             &                                          & $\left( 1-\frac{1}{\rho} , 0 \right) $ \\ \cline{2-4} 
      & \multirow{2}{*}{  $ \hbeta > \mu \rho$ , $\tilde{q}(\htheta, \hpsi) < 1 $ } & $\hbeta > \rho^2 \mu*$  &  $\left (0,   1- \sqrt{\frac{ \mu}{\hbeta} }\right )$    \\ \cline{3-4} 
         &         & $\hbeta < \rho^2 \mu $ & $\left(\frac{\mu \rho}{\hbeta} - \frac{1}{\rho}, 1- \frac{\mu \rho}{\hbeta} \right )$      \\ \cline{2-4}  
         &   
         \multirow{2}{*}{$\hbeta > \rho^2 \mu$, \ 
         ${\tilde q} \left(0,   1- \sqrt{\frac{ \mu}{\hbeta} } \right) > 1$} 
         & $\mu + 1 > \mu \rho$ & $\left(0, \frac{1}{\mu+1} \right)$ \\ \cline{3-4}
         &        & $\mu + 1 < \mu \rho$ & $\left(\theta_E, \psi_E \right )$ \\ \cline{2-4}
 &  \begin{tabular}[c]{@{}c@{}} $\mu \rho < \hbeta < \rho^2 \mu$, $\tilde{q} \left(\frac{\mu \rho}{\hbeta} - \frac{1}{\rho}, 1- \frac{\mu \rho}{\hbeta} \right ) >1$\\
 $\implies$ $\mu \rho > (\mu+1)$ \end{tabular} &  & $\left(\theta_E, \psi_E \right )$ \\ \cline{1-4}
SE, $\rho < 1$ &  $\mu > \hbeta$  & & $(0, 0)$   \\ \hline
\end{tabular}
\vspace{1mm}
\caption{Attractors for FR agents, 
$(\theta_E, \psi_E) := \left(1-\frac{1}{\rho}-\frac{1}{\mu \rho} ,
\ \frac{1} { \mu \rho} \right )$
\label{table_FR}}
    \vspace{-8mm}
\end{table}
As seen from the   two theorems, we have  two types of attractors: (i) interior attractors (for e.g., third row in Tables \ref{table_FC}-\ref{table_FR}) in which $(\htheta, \hpsi) \in (0,1)\times (0, 1)$, and  (ii) boundary attractors, where at least one of the components is $0$. In the latter case, either the disease is eradicated with the help of vaccination ($\htheta = 0$, $\hpsi >0$) or  the disease gets cured without the help of vaccination $(\htheta = \hpsi = 0)$ or no one gets vaccinated $(\hpsi = 0)$. In the last case, the fraction of infected population reaches maximum possible level for the given system, which {\it we refer to as non-vaccinated disease fraction (NVDF), $(\htheta^N = 1 - \nicefrac{1}{\rho})$ (first row in Tables \ref{table_FC}-\ref{table_FR}).} Further, $\hat{\eta}$ is always in $(0,1)$. Furthermore,  important characteristics of the attractors depend upon  quantitative parameters describing the nature and the spread of the disease, and the  vaccination responses:  %
\vspace{-2mm}
\begin{itemize}
    \item \textit{Endemic disease:} The disease is not self-controllable with $\rho > 1$ and the eventual impact of the disease is governed by the  attitude of agents towards vaccination (all rows other than the last in Tables \ref{table_FC}-\ref{table_FR}). 
        \item \textit{Self-eradicating (SE) Disease:} The disease is not highly infectious ($ \rho < 1$) and  can be eradicated without exogenous aid (vaccine). 
\end{itemize}
The above characterisation interestingly {draws parallels from  queuing theory.} One can view $\lambda$ as the arrival of infection and $r+b+d_e$ as its departure. Then  $\rho = \nicefrac{\lambda}{(r+b+d_e)}$,   resembles the load factor. It is well known that queuing systems are stable when $\rho < 1$, similarly in our case, the disease gets self-eradicating with $\rho < 1$.   
\renewcommand{\arraystretch}{1.8} 
 The attractors for VFC1 agents (proof in Appendix B):

\renewcommand{\arraystretch}{1.8} 
\begin{theorem}{\bf [VFC1 agents]}\label{thrm_VFC1}
When    $\hbeta \leq 2\mu\rho^2$ or when $\hat{q} = 1$ with $\tilde{q} (\htheta, \hpsi) \ne  1$,  the attractors for ODE \eqref{eqn_ODE} are $(\htheta, \hpsi, \nicefrac{(b-d)}{{\hat \varrho}})$, and are provided in Table \ref{table_VFC1}. The interior attractors are the zeroes of the RHS of  ODE \eqref{eqn_ODE}. \eop
\end{theorem}

 

\begin{table}[htbp]
\vspace{-5mm}
    \centering
\begin{tabular}{|c|c|c|}
\hline
Nature                     & Parameters                                  & $(\htheta, \hpsi)$  \\ \hline \hline
\multirow{3}{*}{\begin{tabular}[c]{@{}c@{}} Endemic, \\ $\rho > 1$\end{tabular}}  & $\hbeta < \mu \left(\frac{\rho^2}{\rho - 1}\right)$   & $\left(1-\frac{1}{\rho}, 0 \right)$     \\ \cline{2-3}                                            
 & $\hbeta > \mu \left(\frac{\rho^2}{\rho - 1}\right)$,   $\  \  \tilde{q}(\htheta, \hpsi) < 1 $   & $\left(\frac{\mu \rho}{\hbeta}, 1 - \frac{1}{\rho} - \frac{\mu \rho}{\hbeta} \right)$   \\ \cline{2-3} 
    & $\tilde{q} \left(\frac{\mu \rho}{\hbeta}, 1 - \frac{1}{\rho} - \frac{\mu \rho}{\hbeta} \right)   > 1 $ $\implies$ $\rho \mu > \mu + 1$ and $\hbeta > \frac{(\rho \mu)^2}{\rho \mu - \mu - 1}$ & $\left(\theta_E, \psi_E \right )$     \\ \hline
SE, $\rho < 1$                                                         &  &  $(0, 0)$      \\ \hline
\end{tabular}
\caption{Attractors for VFC1
agents: Disease never gets eradicated with $\rho >1$ \label{table_VFC1}}


\vspace{-9mm}
\end{table}

\noindent
{\bf Key observations and comparisons} of the various equilibrium states:\\  
\noindent$\bullet$ When the disease is self-eradicating, agents need  not get vaccinated   to  eradicate the disease. For all the type of agents, $\htheta = 0$ (and so is $\hpsi = 0$) for all $\hbeta < \mu.$

\noindent$\bullet$ With  endemic disease, {\it it is possible to eradicate the disease  only if the agents get vaccinated aggressively}. The FC agents with  $\hbeta > \rho \mu$ and FR agents with   a bigger $\hbeta > \rho^2 \mu$ can completely eradicate the disease; this is possible only when $\mu \rho < \mu + 1$. However  interestingly, \textit{these parameters can't drive the system to an equilibrium that is stable against mutations, as will be seen in the next section.  }
    
\noindent$\bullet$ With a lot more  aggressive FR/FC agents, the system reaches disease free state ($\htheta = 0$ for all bigger $\hbeta$), however  with  bigger vaccinated fractions ($\leq \nicefrac{1}{(\mu+1)}$).
    
\noindent$\bullet$ Interestingly {\it such an eradicating equilibrium state is not observed with vigilant agents (Table \ref{table_VFC1}).}  This is probably analogous to the well known fact that  the rational agents often pay high price of anarchy.  

\begin{figure}[htbp]
\vspace{-35mm}
\begin{minipage}{0.49\linewidth}
    \centering
    \includegraphics[width=6cm,height=7.35cm,keepaspectratio]{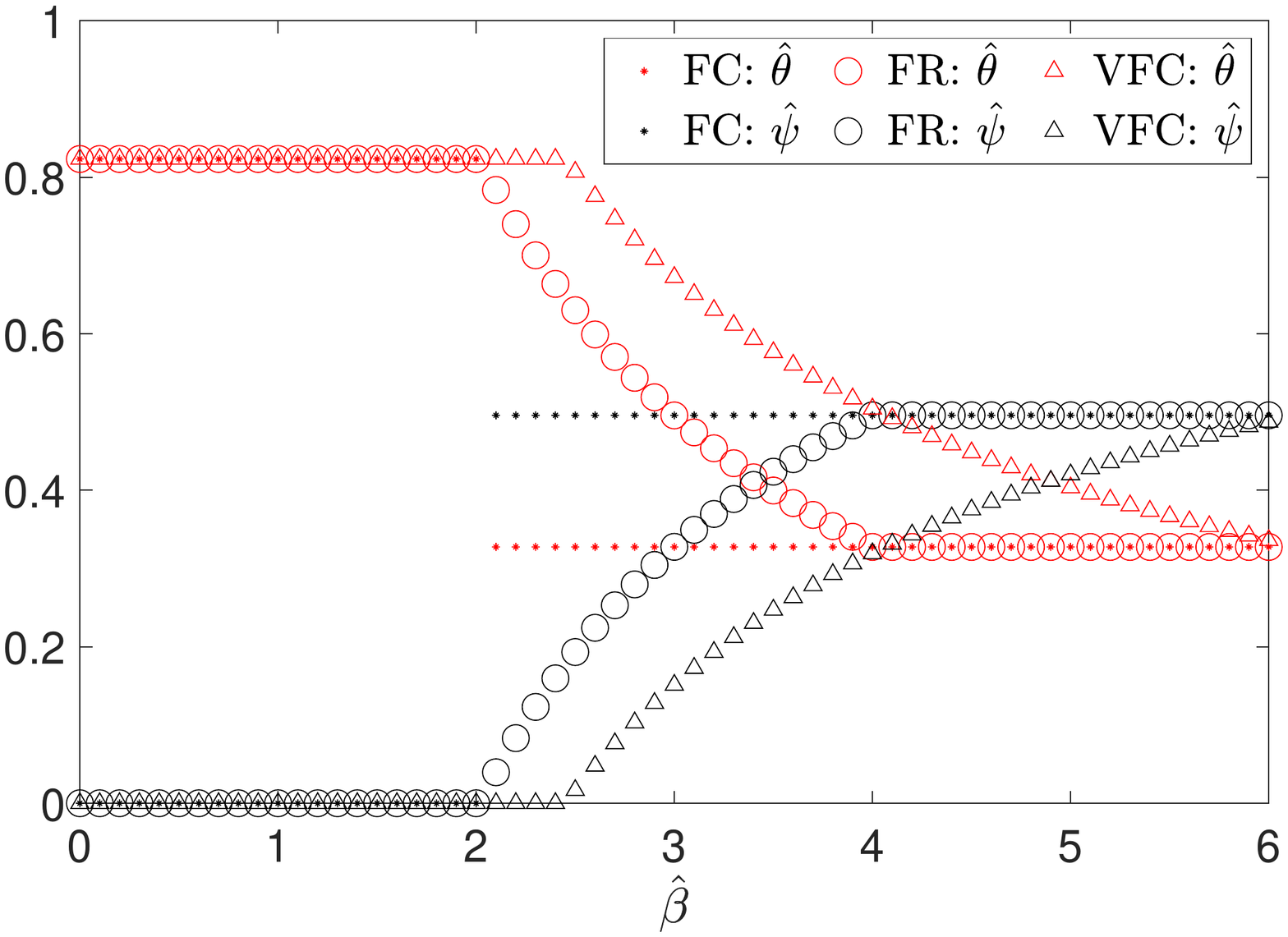}
    \vspace{-35mm}
\end{minipage}%
\begin{minipage}{0.49\linewidth}
    \centering
    \vspace{4mm}
    \includegraphics[width=6.8cm,height=8.9cm]{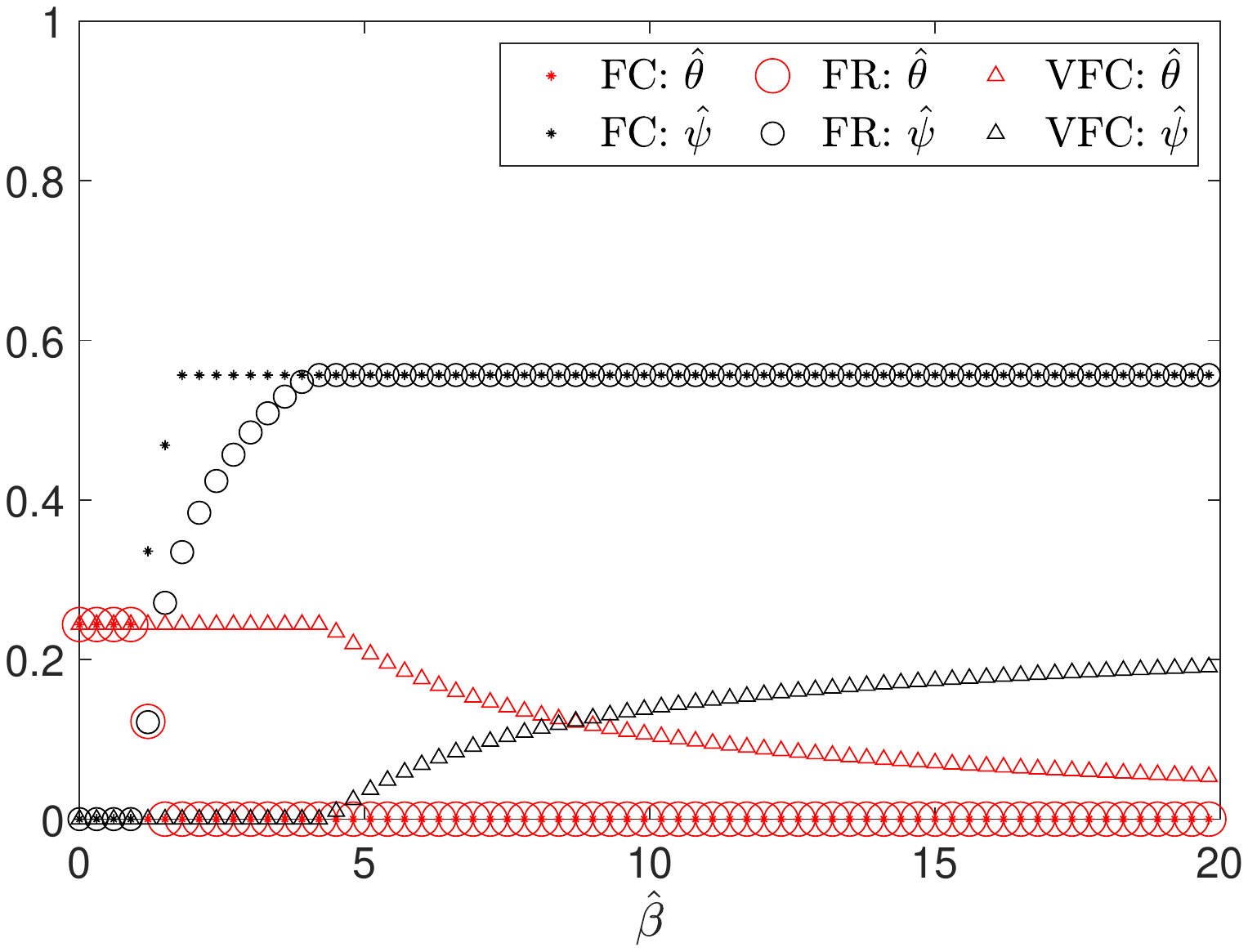}
    \vspace{-35mm}
\end{minipage}
\vspace{2mm}
\caption{Attractors for FC, FR, VFC1 agents versus $\beta$ \label{fig_equlibrium_states}}
\vspace{-8mm}
\end{figure}

\noindent$\bullet$ For certain behavioural parameters, system reaches an equilibrium  state at which vaccinated and infected population co-exist. Interestingly for all three types of agents some of the co-existing equilibrium are  exactly the same  (e.g., $(\theta_E, \psi_E)$ in Tables and left plot of Figure \ref{fig_equlibrium_states}). In fact, such equilibrium  are stable against mutations, as will be seen in the next section. 


A numerical example is presented in Figure \ref{fig_equlibrium_states} that depicts many of the above observations. The parameters in respective plots are ($r =1.188, \nu=0.904$, $ \lambda= 8.549$) and ($r = 1.0002,  \nu=0.404, \lambda= 1.749$), with $b=0.322$. From the left plot, it can be seen that for all agents, $(\htheta, \hpsi)$  equals NVDF for smaller values of $\hbeta$. As $\hbeta$ increases, proportions for FC agents directly reach $(\theta_E, \psi_E)$.  However, with FR and VFC1 agents,
the proportions gradually shift from another interior attractor to finally settle at  $(\theta_E, \psi_E)$. 
Further,  for the right plot ($\mu \rho < \mu + 1$) the disease is eradicated with FC (when $\hbeta > \rho \mu$), FR (when $\hbeta > \rho^2 \mu$), i.e., $(\htheta, \hpsi)$ settles to $(0, \nicefrac{1}{\mu + 1})$, which is not a possibility in VFC1 agents. For the latter type, the proportions traverse through an array of co-existence equilibrium, and would approach $(0, 1-1/\rho)$ as $\hbeta$ increases (see row 2 in Table \ref{table_VFC1}). 



 
\begin{table}[htbp]
\vspace{-5mm}
\centering

 
\scalebox{0.9}{
\begin{tabular}{|c|c|c|c|}
\hline
Nature                    & \multicolumn{2}{c|}{Parameters}                                                                                                                 & $(\htheta, \hpsi)$ \\ \hline \hline
\multirow{3}{*}{\begin{tabular}[c]{@{}c@{}} Endemic, \\ $\rho > 1$\end{tabular}}  & $ \hbeta > \rho \mu$   *   &                                                      & $\left( 0, 1-\frac{\mu}{\hbeta} \right) $                                                                                \\ \cline{2-4} 
                                       & \multirow{3}{*}{$ \hbeta < \rho \mu$} & $ \rho_e \mu_e > 1$  or $\hbeta \nu < b- d_e$ & $\left( 1-\frac{1}{\rho_e} , 0\right)$                                                                   \\ \cline{3-4} 
                                       &                                                                                          &$ \rho_e \mu_e < 1$,  $\hbeta \nu > b -d_e$ & \begin{tabular}[c]{@{}c@{}}$\left (\theta^*,   1- \theta^*\left(1-\frac{ d_e}{\lambda}\right) -\frac{1}{\rho}\right )$,   \\  $\theta^* = \frac{\hbeta \nu \left(\frac{\mu}{\hbeta} -  \frac{1}{\rho_e}\right)}{d_e \left(1-\frac{\hbeta \nu}{\lambda} \right)}$ \end{tabular}                     \\ \hline
SE, $\rho < 1$                      & $\hbeta < \mu$                                      &                                                      & $(0, 0)$                                                                                                                   \\ \hline
\end{tabular}}
\caption{Attractors for  FC agents (deadly disease) \label{table_FC_deadly}}
\vspace{-8mm}
\end{table}

\noindent {\bf Deadly disease:}\label{sec_deadly_disease} 
We now consider the deadly disease scenario ($d_e > 0$) with FC and FR agents. 
Let {\small
$
\rho_e := \frac{\lambda - d_e}{r+b},  \mu_e := \frac{b - d_e}{\hbeta \nu - d_e}.
$}
We  conjecture 
the attractors  with  $\hbeta \le 1$ in Tables \ref{table_FC_deadly}, \ref{table_FR_deadly} respectively,  with $A := d_e \hbeta \nu,\ B := -[(r+b+d_e)\hbeta \nu + d_e \lambda], \mbox{ and }\  C := b \lambda + d_e (r + d_e). $
Further, the candidate attractors with   $\tilde{q} (\htheta, \hbeta) > 1$ are provided in the next section, which is of interest to ESS-AS. We expect that the proofs can be extended analogously and will be a part of future work, along with identifying and proving the attractors for other cases.

\begin{table}[htbp]
\vspace{-3mm}
\centering

 
\scalebox{0.9}{
\begin{tabular}{|c|c|c|c|}
\hline
Nature                  & \multicolumn{2}{c|}{Parameters}                                                                                                                 & $(\htheta, \hpsi)$ \\ \hline \hline
\multirow{6}{*}{{\begin{tabular}[c]{@{}c@{}} Endemic, \\ $\rho > 1$\end{tabular}} } &  \multirow{3}{*}{$ \hbeta > \rho \mu$ } & $\hbeta > \rho^2 \mu \ *$  & $\left( 0, 1 - \sqrt{\frac{\mu}{\hbeta}}\right)$                                                                   \\ \cline{3-4} 
                                       &                                                                                          & $\hbeta < \rho^2 \mu$ & \begin{tabular}[c]{@{}c@{}}$\left (1 - \frac{1}{\rho_e} - \frac{\lambda \hpsi}{\lambda-d_e}, \hpsi \right )$,   \\  $\hpsi = 1 - \frac{-B - \sqrt{B^2 - 4AC}}{2A} $ \end{tabular}                     \\ \cline{2-4} 
                                       & \multirow{3}{*}{$ \hbeta < \rho \mu$} & $ \rho_e \mu_e > 1$  or $\hbeta \nu < b-d_e$ & $\left( 1-\frac{1}{\rho_e} , 0\right)$                                                                   \\ \cline{3-4} 
                                       &                                                                                          &$ \rho_e \mu_e < 1$,  $\hbeta \nu >b-d_e$ & \begin{tabular}[c]{@{}c@{}}$\left (1 - \frac{1}{\rho_e} - \frac{\lambda \hpsi}{\lambda-d_e}, \hpsi \right )$,   \\  $\hpsi = 1 - \frac{-B - \sqrt{B^2 - 4AC}}{2A} $ \end{tabular}                    \\ \hline
SE, $\rho < 1$                      & $\hbeta < \mu$                                      &                                                      & $(0, 0)$                                                \\ \hline
\end{tabular}}
\caption{Attractors for FR agents (deadly disease) \label{table_FR_deadly}}
\vspace{-14mm}
\end{table}

\section{Evolutionary stable vaccination responses}
Previously, for a given user behaviour, we showed that the system reaches an equilibrium state, and identified the corresponding equilibrium states in terms of limit proportions. If such a system is invaded by mutants that use a different vaccination response, the system can get perturbed, and there is a possibility that the  system drifts away. We now identify those equilibrium states, which are evolutionary stable against static mutations. Using standard tools of evolutionary game theory, we will show that the mutants do not benefit from deviating under certain subset of policies. That is, we identify the ESS-AS policies defined at the end of  Section \ref{sec_prob_desc}.

We begin our analysis with $d_e=0$, the case with no excess deaths (due to disease).  
From Tables \ref{table_FC}-\ref{table_VFC1},  with $\rho > 1$,  for all small values of $\hbeta$ (including $\hbeta=0$), irrespective of the type of the policy, the equilibrium state remains the same (at NVDF), $(\htheta,\hpsi)=\left(1-\nicefrac{1}{\rho}, 0\right)$. Thus, the  value of $h(\cdot)$ in user utility function \eqref{eqn_utility} for all such small $\hbeta$ equals the same value and is given by:

\vspace{-5mm}
{\small 
\begin{eqnarray}\label{eqn_cost}
    h_m := h(\pi(0)) =h\left(1-\frac{1}{\rho}, 0\right)=  c_{v_1} +\bar{c}_{v_2} - p_I\left(\theta \right) c_{I_1}, 
    \mbox{ with, }     p_I\left(\theta\right) := \frac{\lambda \theta }{\lambda\theta +\nu}.\hspace{4mm} 
\end{eqnarray}}In the above, $p_I(\cdot)$ is the probability that the individual gets infected before the next vaccination epoch.  
The  quantity $h_m$  is instrumental in deriving the following result with $\rho>1$. When $\rho < 1$, the equilibrium state for all the policies (and all $\hbeta$) is $(0,0)$ leading to the following (proof in \cite{TR}):
\begin{lemma}\label{lem_ess_zero}
If $\rho<1$, or if $\rho  > 1$ with $h_m>0$, then $\pi({0})$ is an ESS-AS, for any  $\pi \in \Pi$. \eop
\end{lemma} 
When the disease is self-eradicating ($\rho<1$), the system converges to $(0,0)$, an infection free state on it's own without the aid of vaccination. Thus we have the above ESS-AS. When $\rho>1$, if the inconvenience caused by the disease  captured by $-h_m $ is not compelling enough (as $-h_m<0$), the ES equilibrium state again results at $\hbeta=0$.    \textit{In other words, policy to never vaccinate is evolutionary stable in both the cases}. Observe this  is a static policy irrespective of agent behaviour (i.e., for any $\pi \in \Pi$), as with $\hbeta= 0$ the agents never get vaccinated irrespective of the system state. 

From  tables of  the previous section, there exists $\bar{\beta}(\pi)$ such that the equilibrium state remains at NVDF for all $\hbeta<\bar{\beta} $  (including  $\hbeta=0$) and for all $\pi \in \Pi$. For   such $\hbeta$, we have $h(\htheta,\hpsi)<0$   if $h_m<0$. Thus from the user utility function \eqref{eqn_utility} and ESS-AS definition, the static best response set ${\cal B} (\hpi)=\{1\}$ for all $\pi\in\Pi$ and all $\hbeta<\bar{\beta}$. Thus with $\rho>1$ and $h_m<0$, any policy $\hpi$ such that $\hbeta<\bar{\beta}$ is not an ESS-AS. This leads to the following   \TR{(proof  available in \cite{TR})}{}:

\begin{theorem}\label{Thrm_ESS}
{\bf[Vaccinating-ESS-AS]}
When  $\rho >1$ and $h_m < 0$, there exists an ESS-AS among a  $\pi\in\Pi$ if and only if
the following two conditions hold:
\begin{enumerate}[(i)]
    \item there exist a $\beta^*>0$ such that $q(\up^*)=1$ and $\tilde{q}( \up^*)\ne 1$ under policy $\pi(\beta^*)$,
    \item   the equilibrium state is
    {\small$
(\theta^*,\psi^*)=\left(\theta_E,\psi_E \right )\hspace{-1mm},$} with $\mu \rho>\mu+1$ and the corresponding user utility component, 
    {\small $h(\pi(\beta^*))= h\left(\theta_E,\psi_E \right )<0$}. \eop
\end{enumerate}

\end{theorem}
\TR{}{
\textbf{Proof:}  Say $\hpi$ is an ESS-AS, for some $\hbeta>\bar{\beta}$. For $\rho>1$, and $\hbeta>\bar{\beta}>0$, from Tables \ref{table_FC}-\ref{table_VFC1},  $q^*:=q(\htheta,\hpsi)\ne 0$. Also,
$\{q^*\}={\cal B}(\hpi)$, by definition of ESS-AS. Since the   utility  \eqref{eqn_utility}   is linear in $q$,   there are only two possible minimizers, $0$ or $1$ depending on $h(\hbeta)$ at ESS (because the BR set at ESS has unique element). So, if $\hpi$ is ESS-AS, then $q^*=1$, 
which implies $h(\hpi)<0$. From the ODE \eqref{eqn_ODE}  with $q^* = 1$,  it is clear that the zero of the RHS of ODE and hence the  attractor (using Theorems \ref{thrm_FCagents}-\ref{thrm_VFC1}) is $\left(\theta_E,\psi_E \right )$. Observe $\theta_E>0$ if and only if $\mu\rho>\mu+1$. 

 Given the possibility of a  $q^*=1$ and $h(\hpi) = h\left(  \theta_E,\psi_E \right ) <0$, first requirement of ESS-AS is satisfied with $\{q^*\}={\cal B}(\hpi)$. Further, mutational stability follows from Lemma \ref{lem_cont_of_attract} because $\tilde{q}(\up^*)\ne 1$.
\eop}

{\bf Remarks:}
Thus when $\rho>1$ and $h_m<0$, there is no  ESS-AS for any $\pi\in\Pi$ if $h(\theta_E,\psi_E)\ge 0$; observe  that $(\theta_E,\psi_E)$ is equilibrium state with $\hat{q}=1$ and hence from ODE \eqref{eqn_ODE},  does not depend upon $\pi$. On the contrary, if $h(\theta_E,\psi_E)<0$, $\pi(\beta^*)$ is ESS-AS  for any $\pi \in \Pi$, with  $\beta^*>\mu \rho $, $\beta^*>\nicefrac{(\mu\rho)^2}{(\mu \rho -1)}$, and
$\beta^*>\nicefrac{(\mu\rho)^2}{(\mu \rho -1-\mu)}$ respectively for $FC$, $FR$ and $VFC$ policies (using Tables \ref{table_FC}-\ref{table_VFC1}). 

Thus interestingly, evolutionary stable behaviour is either possible in all, or in none. However, the three types of dynamic agents require different set of parameters to arrive  at ES equilibrium. \textit{An ES equilibrium with vaccination is possible only when $\mu\rho> \mu +1$} and interestingly, the infected and vaccinated fractions at this equilibrium are indifferent of the agent's behaviour.
\textit{In conclusion the initial dynamics could be different under the three different agent behaviours, however,  the limiting proportions corresponding to any ESS-AS are the same.}

\subsubsection{Numerical examples:} 
\begin{figure}
\vspace{-12mm}
\begin{minipage}{5cm}
      \centering
    \includegraphics[width=6cm, height = 5.5cm]{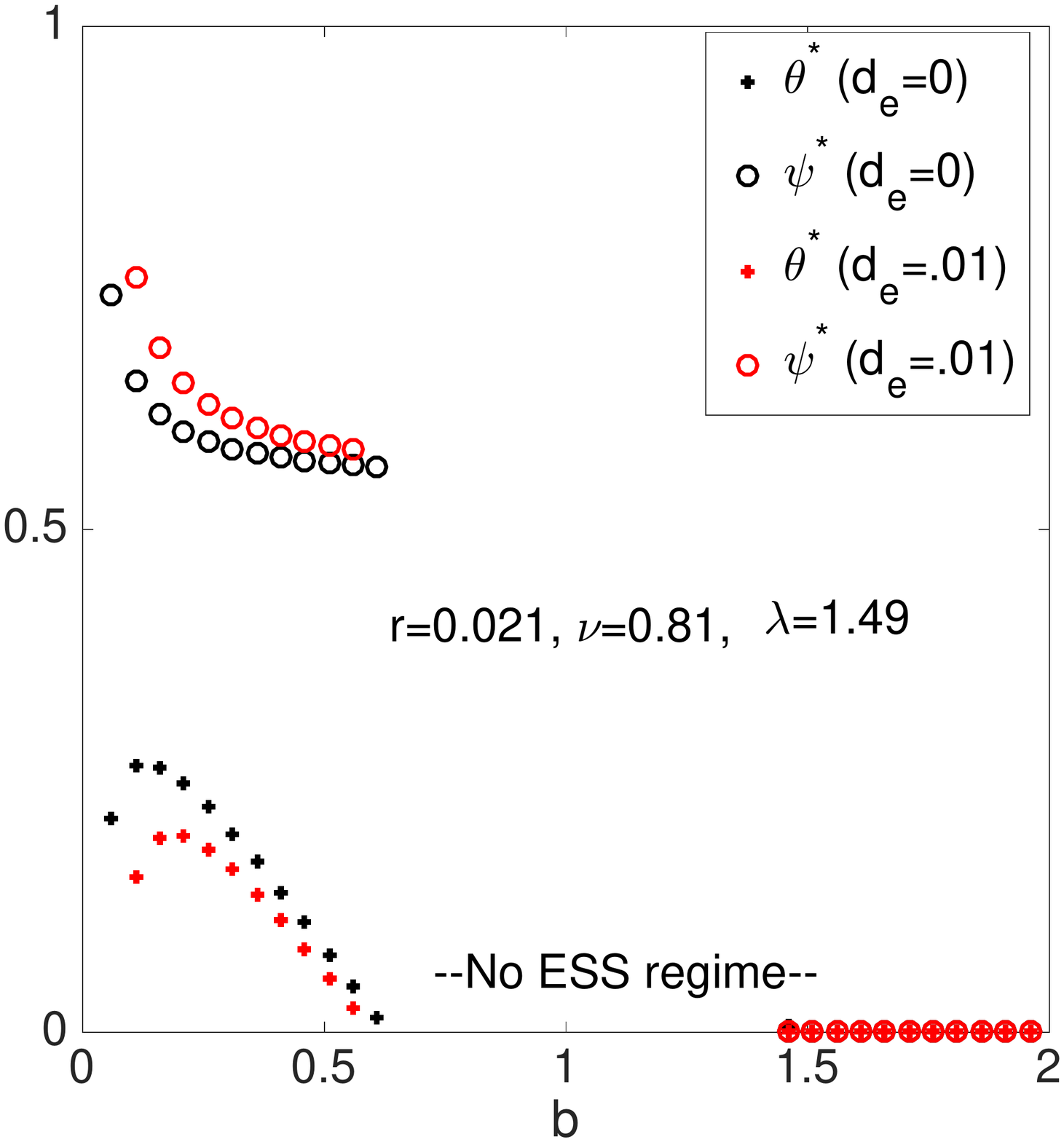}
    \vspace{-17mm}
    \caption{ESS versus birth-rate}
    \label{fig_ess_vs_b}  
\end{minipage}
\hspace{10mm}
\begin{minipage}{5.5cm}
      \centering
    \includegraphics[width=5
    cm, height = 5.5cm]{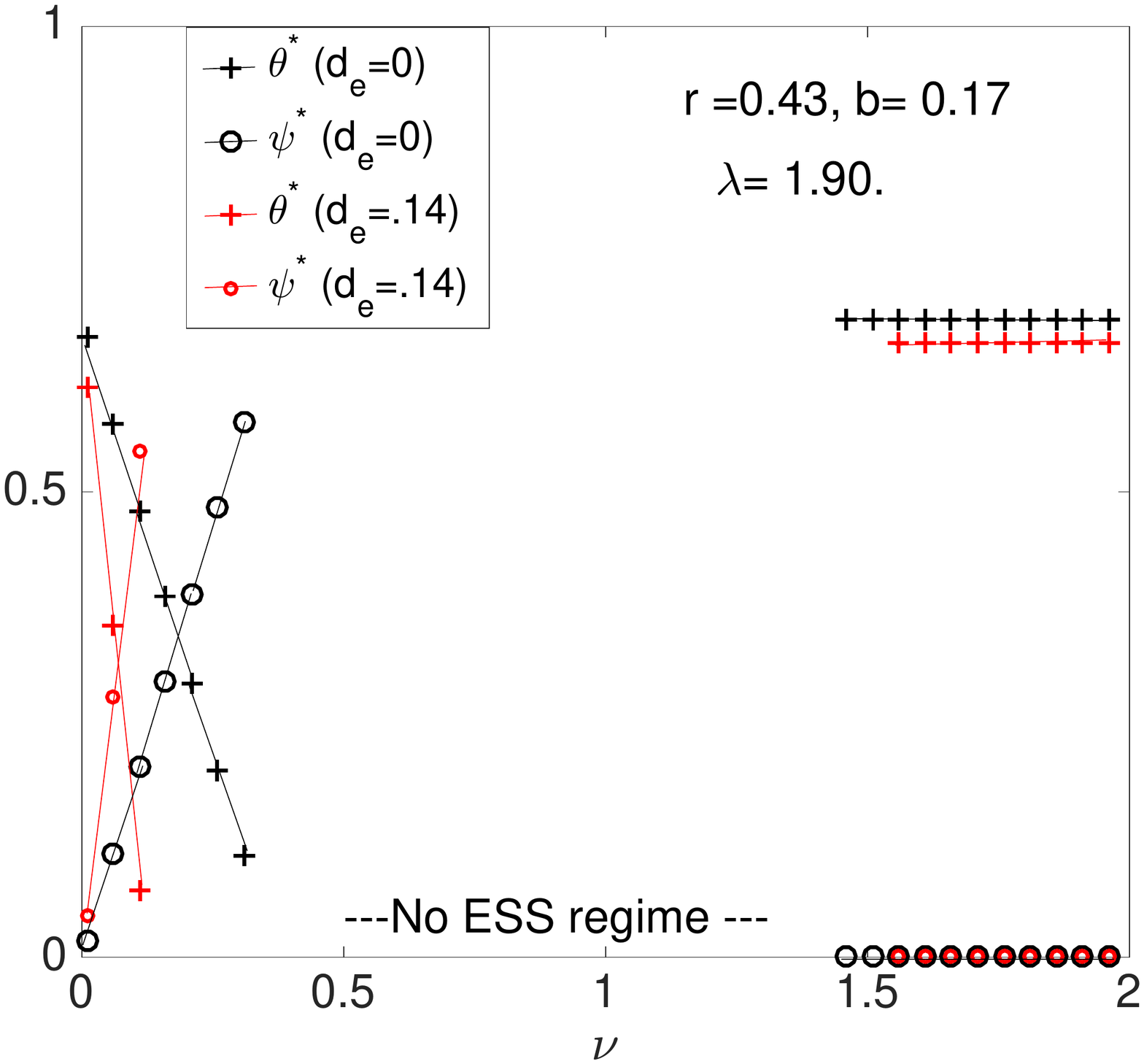}
      \vspace{-12mm}
    \caption{ESS versus vaccine availability}
    \label{fig_ess_vs_nu}  
\end{minipage}
\TR{\vspace{-7mm}}{\vspace{-4mm}}
\end{figure}

We study the variations in vaccinating ESS $(\theta_E,\psi_E)$, along with others, with respect to different parameters. In these examples, we set the costs of vaccination and infection as $c_{v_1}=2.88, c_{v_2}=0.65, \bar{c}_{v_2}=1.91, c_{I_1}=4.32/r$. Other parameters are in the respective figures; black curves are for $d_e=0$.
In Figure \ref{fig_ess_vs_b}, we plot the  ESS-AS   for different values of birth-rate. Initially, $\rho>1,$ and vaccinating ESS-AS exists
for all $b \le 0.54$; here   $h(\theta_E,\psi_E)<0$ as given by Theorem \ref{Thrm_ESS}. As seen from the plot, $\theta_E$ is decreasing and approaches zero at $b\approx 0.54$. Beyond this point there is no ESS because  $\mu\rho$ reduces below $ \mu+1$. With further increase in $b$, non-vaccinating ESS emerges as $\rho$ becomes less than one. Interestingly, a much larger fraction of people get vaccinated at ESS for smaller birth-rates. This probably could be  because of higher infection rate per birth. In fact from the definition of $\theta_E$, the infected fractions at ES equilibrium are concave functions of birth rate. When infection rate per birth is sufficiently high, {\it it appears people pro-actively vaccinate themselves, and bring infected fraction (at ES equilibrium) lower than those   at smaller ratios of infection rate per birth.}

In Figure \ref{fig_ess_vs_nu} we plot ESS-AS for different values of $\nu$. 
For  all $\nu\le 0.31$, the vaccinating ESS-AS exists.  Beyond this, there is no ESS because $\mu\rho$ reduces below $\mu+1$. As $\nu$ further increases, $h_m$ 
becomes\footnote{This is because the chances of infection before the next vaccination epoch decrease with increase in the availability rate $\nu$.} positive, leading to $NVDF$ as ESS.  One would expect a smaller infected proportion at ES equilibrium with increased availability rate, however we observe the converse; this is because the users' perception about infection cost changes with abundance of vaccines.   


\TR{}{
In Figure \ref{fig_ess_vs_lambda}, we plot ESS for different values of $\lambda$. Initially, $\rho<1$ so $(0,0)$ is the ESS. When $\rho$ increases beyond one but $\mu \rho$ is less than $\mu+1$, there is no ESS-AS. With further increase in $\rho$, $\mu \rho$ becomes more than $ \mu+1$, and we get vaccinating ESS as given by Theorem \ref{Thrm_ESS}. When infection rate is very large as compared to $\nu$, the vaccinated fraction even at vaccinated ESS drops below the infected fraction.

In Figure \ref{fig_ess_vs_r}, we plot ESS for different values of  $r$. Initially $\rho>1$, and in this example, $\mu \rho> \mu+1$, so we have vaccinating ESS. With increase in $r$, the infected fraction decreases, and $h(\htheta,\hpsi)$ becomes non-negative, so there is no ESS. As $r$ further increases, as long as $\rho<1$,  $NVDF$ becomes ES equilibrium. Once $\rho$ drops below 1, disease becomes self-eradicating. As seen from the plot, when recovery rate is low, people vacciante them pro-actively and infections come down. Once the infections start dropping, people estimate the cost of infection to be low, and do not consider vaccination to be the best response. Also, with high recovery rate, disease becomes less severe. This leads to an increase in the infected fraction, but with further increase in $r$, disease becomes self eradicating.

\begin{figure}
\vspace{-7mm}
\begin{minipage}{5cm}
      \centering
    \hspace{2mm}\includegraphics[width=5cm, height = 6cm]{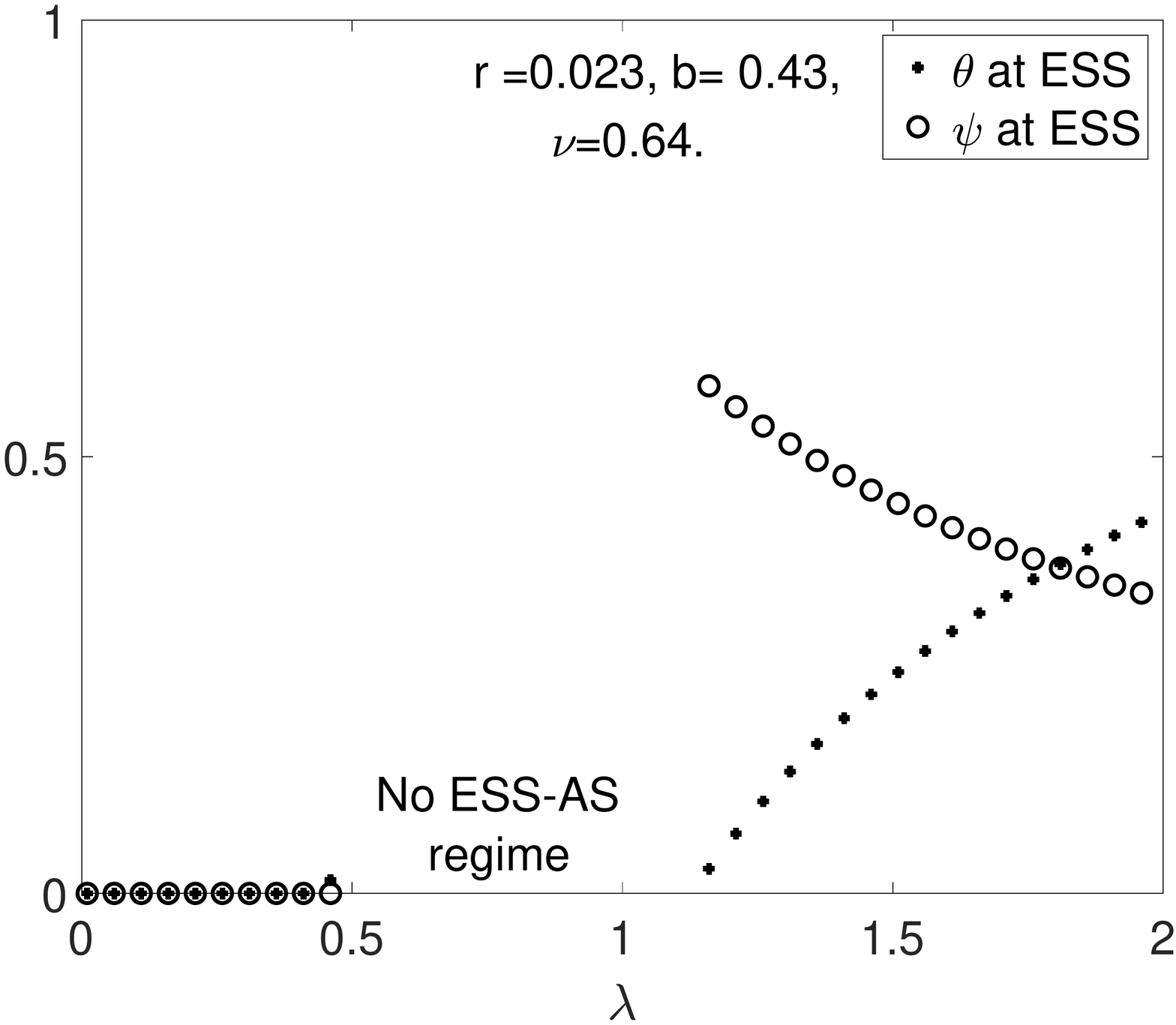}
    \vspace{-15mm}
    \caption{ESS versus infection-rate}
    \label{fig_ess_vs_lambda}  
\end{minipage}
\hspace{10mm}
\begin{minipage}{5cm}
      \centering
      \includegraphics[width=5cm, height = 6cm]{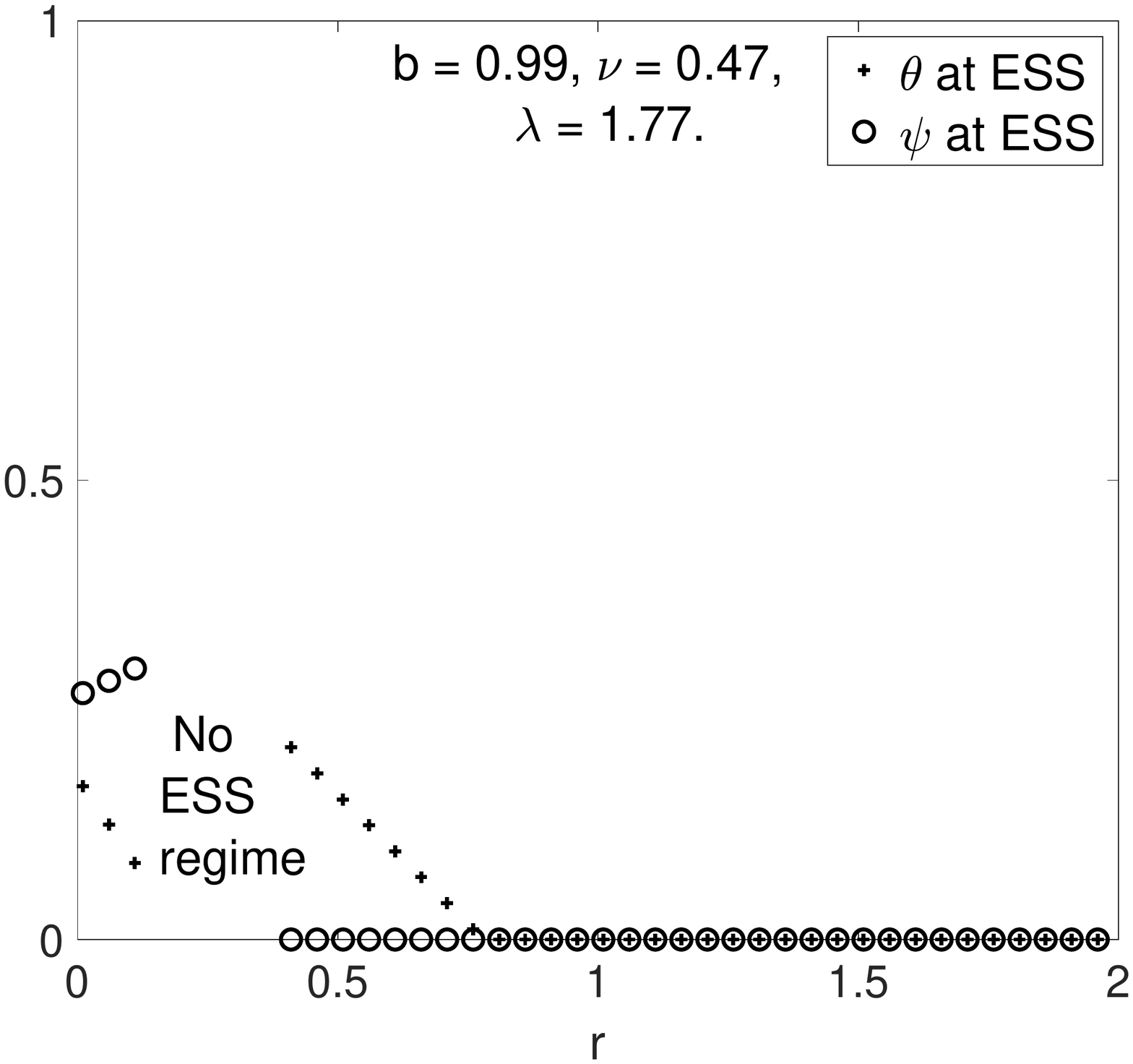}
      \vspace{-15mm}
    \caption{ESS versus recovery rate}
    \label{fig_ess_vs_r}  
\end{minipage}
\vspace{-4mm}
\end{figure}}
\vspace{-3mm}
\subsubsection{With excess deaths ($d_e>0$):} The analysis  will follow in exactly similar lines as above. In this case, we have identified the equilibrium points of the ODE \eqref{eqn_ODE} but are yet to prove that they are indeed attractors. We are hoping that proof of attractors will go through similar to case with $d_e=0$, but we have omitted it due to lack of time and space. Once we prove that the equilibrium points  are attractors, the analysis of ESS would be similar to the previous case. When $\rho<1$  or if $\rho>1$ and  $h_{m}^{d_e}:=h(\pi(0))=h\left(1-\frac{1}{\rho_e},0\right)>0$ (see \eqref{eqn_utility}, Table \ref{table_FC_deadly}-\ref{table_FR_deadly}), then $\pi(0)$ is an ESS-AS, as in Lemma \ref{lem_ess_zero}. Now, we are only left with case when $\rho>1$ and $h_m^{d_e}<0$ and one can proceed as in Theorem \ref{Thrm_ESS}. In this case the only  candidate for ESS-AS is $\hpi$ with $\hbeta$ such that $q=1$ and $\tilde{q}\ne 1$ and $h(\htheta,\hpsi)<0$. So, we will only compute the equilibrium points for  $\hbeta$ such that $q=1$ and $\tilde{q}\ne 1$. From ODE \eqref{eqn_ODE}, such an equilibrium point is given by ($B:=\lambda b+d_e(r+d_e-\lambda-\nu)$): 
{\small
$$ \theta_{E}^{d_e} =1-\frac{1}{ \rho_e}-\frac{\lambda {\psi}_E^{d_e}}{\lambda-d_e} \mbox{ and } {\psi}_E^{d_e}=\frac{-B+\sqrt{B^2+4\lambda d_e \nu (r+b)}}{2\lambda d_e}.$$ }
One can approximate this root for small   $d_e$ (by neglecting second order term $\lambda d_e \approx 0$), 
 the corresponding  ES equilibrium state $(\theta_{E}^{d_e},\psi_{E}^{d_e})$ (again from \eqref{eqn_ODE}):

\vspace{-4mm}
{\small 
\begin{eqnarray*}
\psi_{E}^{d_e} &\approx&
\frac{(r+b)\nu}{\lambda b+d_e(r+d_e-\lambda-\nu)}, \mbox{ that is, } \\
(\theta_{E}^{d_e},\psi_{E}^{d_e}) &\approx& \left ( 1 - \frac{1}{\rho_e} - \frac{o^{d_e}} {\mu \rho_e}, \  
\frac{o^{d_e}}{\mu \rho_e} \frac{\lambda-d_e}{\lambda} \right ) \mbox{ with }
  \ o^{d_e} := \frac{1}{ 1 +\frac{ d_e(r+d_e-\lambda-\nu)}{ \mu\lambda \nu }}
.\end{eqnarray*}}
As before, there  is no ESS if $
\mu + o^{d_e}  \ge  \mu 
\rho_e$ (for larger $d_e$, when $\theta_E^{d_e} < 0$). 
From Figures \ref{fig_ess_vs_b}-\ref{fig_ess_vs_nu} (red curves), 
 the ES equilibrium state in deadly case has higher vaccinated fraction and lower infected fraction as compared to the corresponding non-deadly case (all parameters same, except for $d_e$). More interestingly the variations with respect to the other parameters remain the same as before.

\begin{figure}
\vspace{-9mm}
\begin{minipage}{7.5cm}
     
\begin{minipage}{3.7cm}
      \centering
    \includegraphics[width=3.8cm, height = 4cm]{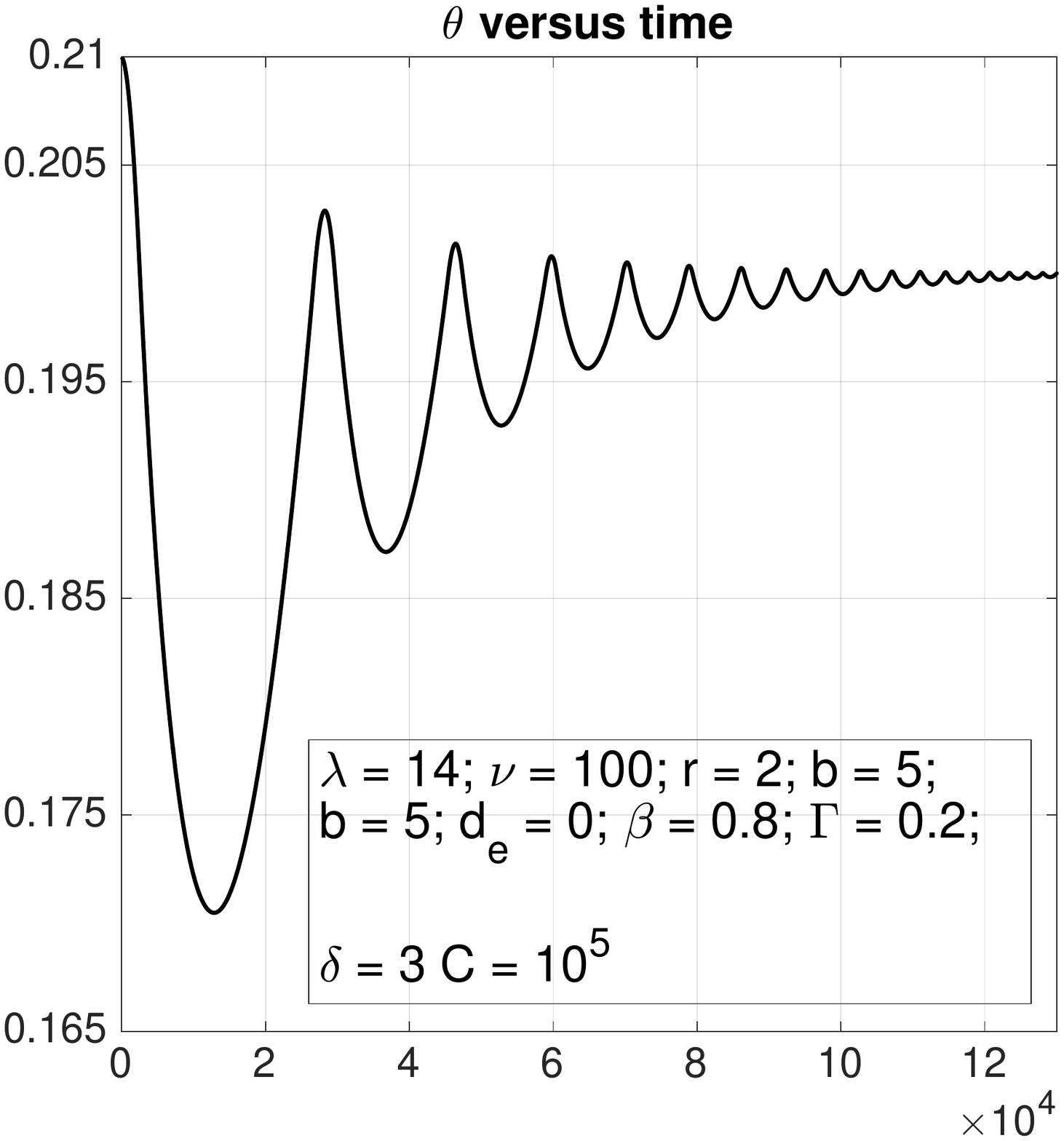}
    \vspace{-10mm}
\end{minipage}
\hspace{.5mm}
\begin{minipage}{3.cm}
      \centering
      \vspace{-1mm}
      \includegraphics[width=3.8cm, height = 3.4cm]{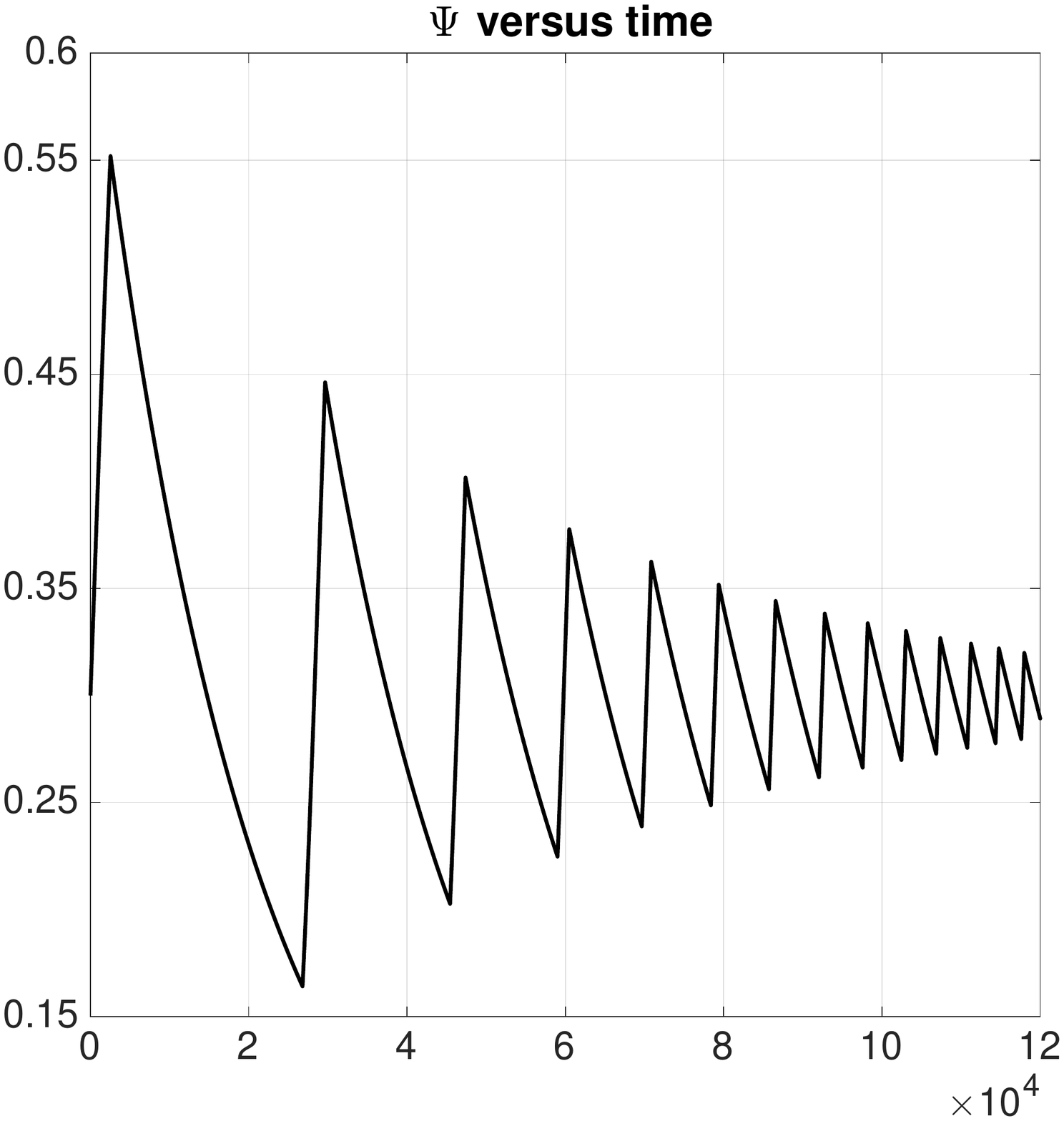}
      \vspace{-8mm}
\end{minipage}
\vspace{-4mm}
 \caption{VFC2 agents: Limit behaviour}
    \label{fig_psi}  
\end{minipage}
\hspace{5mm}
\begin{minipage}{4cm}
      \centering
      \vspace{-12mm}
    \includegraphics[width=3.8cm, height = 5.7cm]{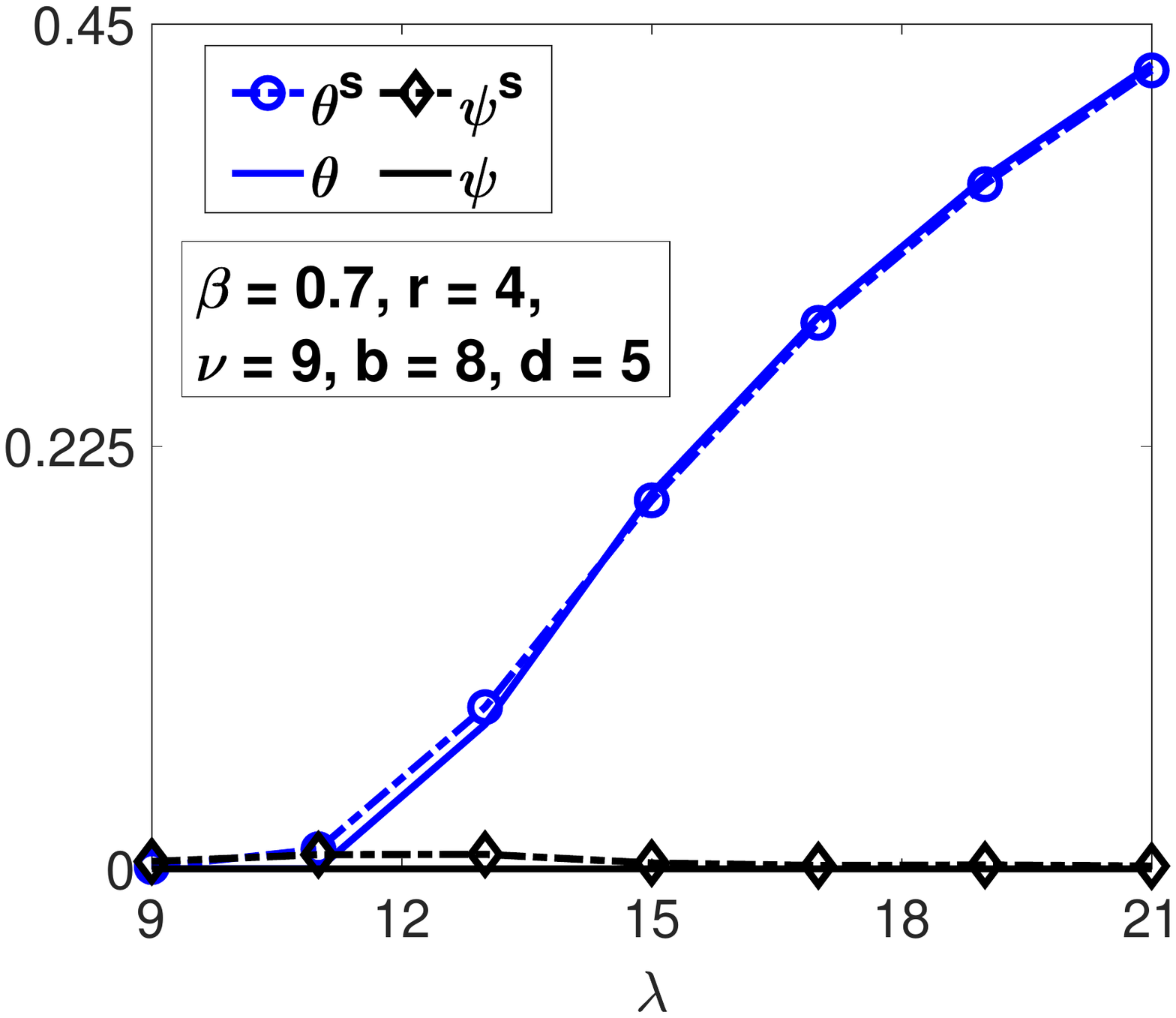}
    \vspace{-18mm}
    \caption{FC agents, against $\lambda$}
    \label{fig_FClambda}  
\end{minipage}
\vspace{-8mm}
\end{figure}

\vspace{-2mm}
\section{Numerical Experiments}

\begin{figure}
\vspace{-18mm}
\begin{minipage}{3.5cm}
      \centering
    \includegraphics[width=3.7cm, height = 6.3cm]{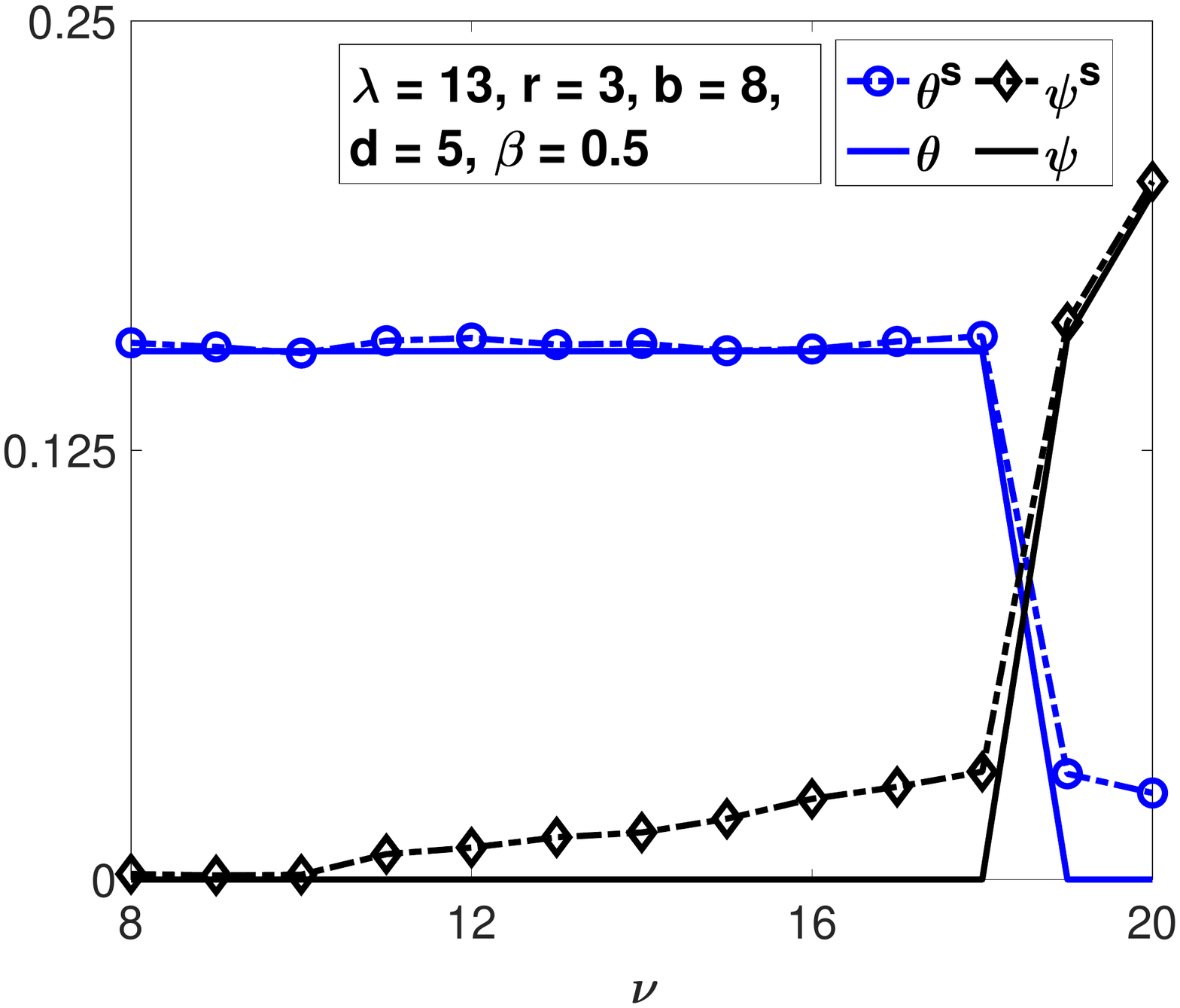}
      \vspace{-18mm}
    \caption{FC agents vs $\nu$}
    \label{fig_FCnu}  
\end{minipage}
\hspace{5mm}
\begin{minipage}{3.5cm}
      \centering
    \includegraphics[width=4cm, height = 6.5cm]{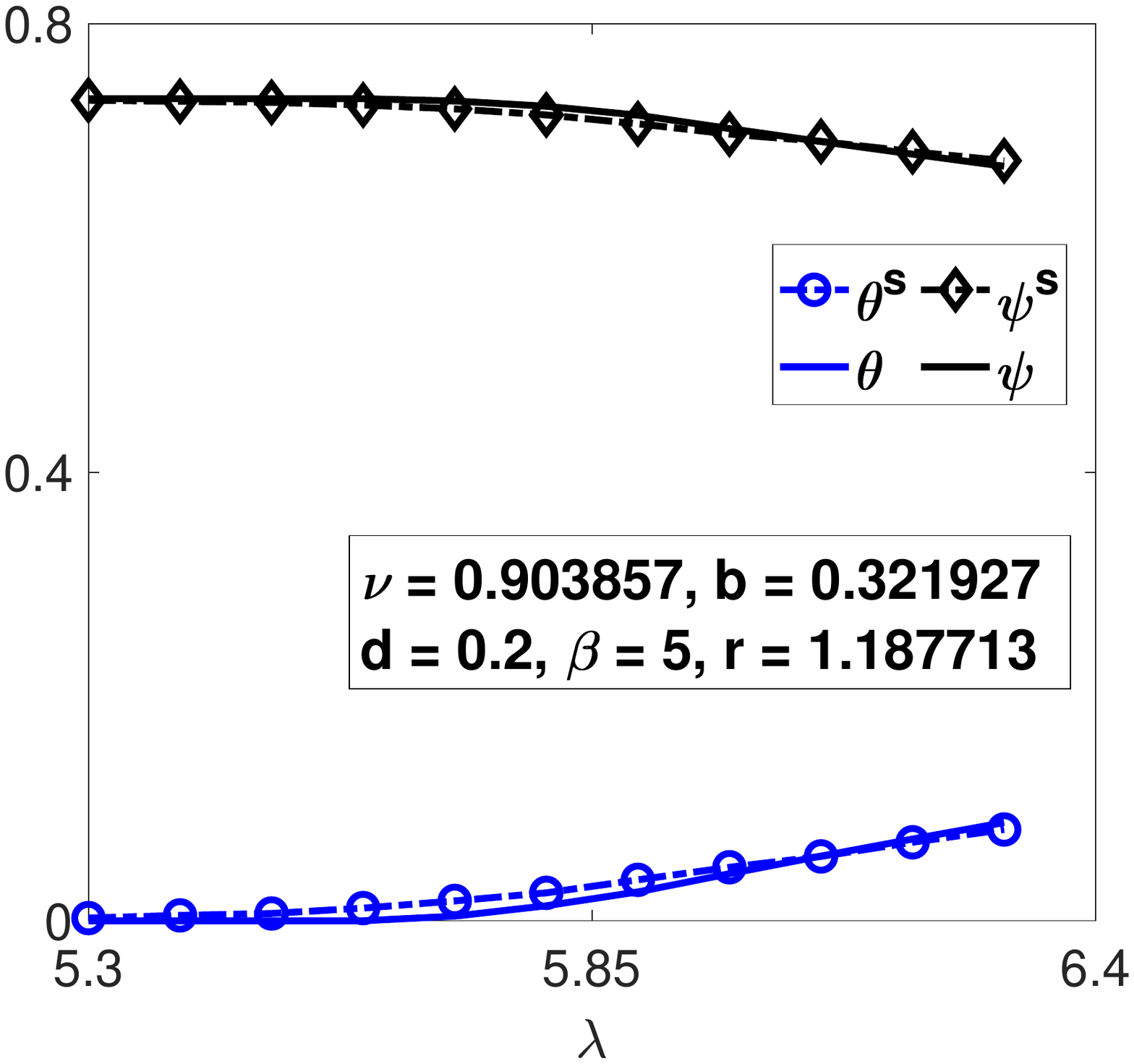}
    \vspace{-20mm}
    \caption{FR agents vs $\lambda$}
    \label{fig_FRlambda}  
\end{minipage}
\hspace{2mm}
\begin{minipage}{3.5cm}
      \centering
     \vspace{-5mm} 
    \includegraphics[width=4.6cm, height = 7.3cm]{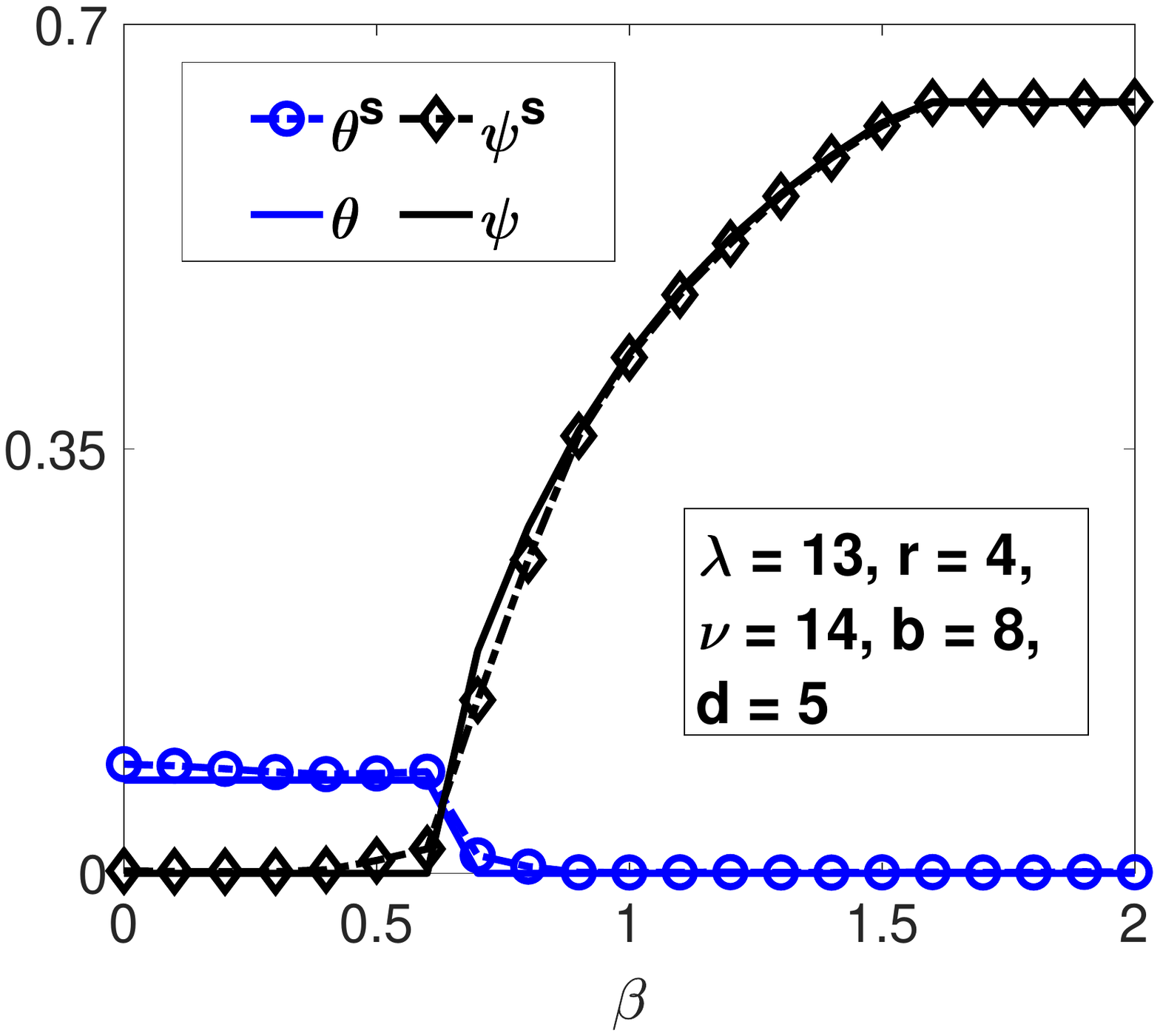}
      \vspace{-23mm}
    \caption{FC agents vs $\beta$}
    \label{fig_FCbeta}  
\end{minipage}
\vspace{-12mm}
\end{figure}

We performed Monte-Carlo simulations to reinforce our ODE approximation theory. We plotted attractors of the ODE \eqref{eqn_ODE} represented by $(\theta,\psi)$, and the corresponding infected and vaccinated fractions $(\theta^s,\psi^s)$ obtained via simulations for different values of $\lambda$, $\nu$ and $\beta$.
Our Monte-Carlo simulation based dynamics mimic the    model described in Section \ref{ODE_approx}.
In all these examples we set $N(0) = 40000$. The remaining parameters are described in the respective figures. We have several plots in figures \ref{fig_FClambda}-\ref{fig_FCbeta}, which illustrate that the  ODE attractors well approximate the system limits, for different sets of parameters.

\noindent{\bf VFC2 agents:}  These agents attempt to vaccinate themselves only when the disease is above a certain threshold $\Gamma$, basically $\tilde {q} = \hat{\beta} \psi \indc{\theta > \Gamma }$. As one may anticipate, the behaviour of such agents is drastically different from the other type of agents. Theorem \ref{thrm1} is applicable even for these agents (approximation in finite windows will be required here). 
However with a close glance at the ODE, one can identify that the ODE does not have a limit point or attractor, but rather would have a limiting set. From the RHS of the ODE \eqref{eqn_ODE},  one can observe that the $\psi$ derivative fluctuates between positive and negative values, and hence $\psi$ goes through increase-decrease phases if at-all  $\theta (t)$ reaches near $\Gamma$. This indeed happens, the fact  is supported by a numerical example of Figure~\ref{fig_psi}.   
Thus interestingly with such a vaccine response behaviour, the individuals begin to vaccinate the moment the infection is above $\Gamma$, which leads to a reduced infection, and when it reaches below $\Gamma$, individuals stop vaccinating themselves. This continues forever, and one can observe such behaviour even in real world.

\hide{
\begin{figure}
\begin{minipage}{5cm}
      \centering
    \includegraphics[width=6cm, height = 8cm]{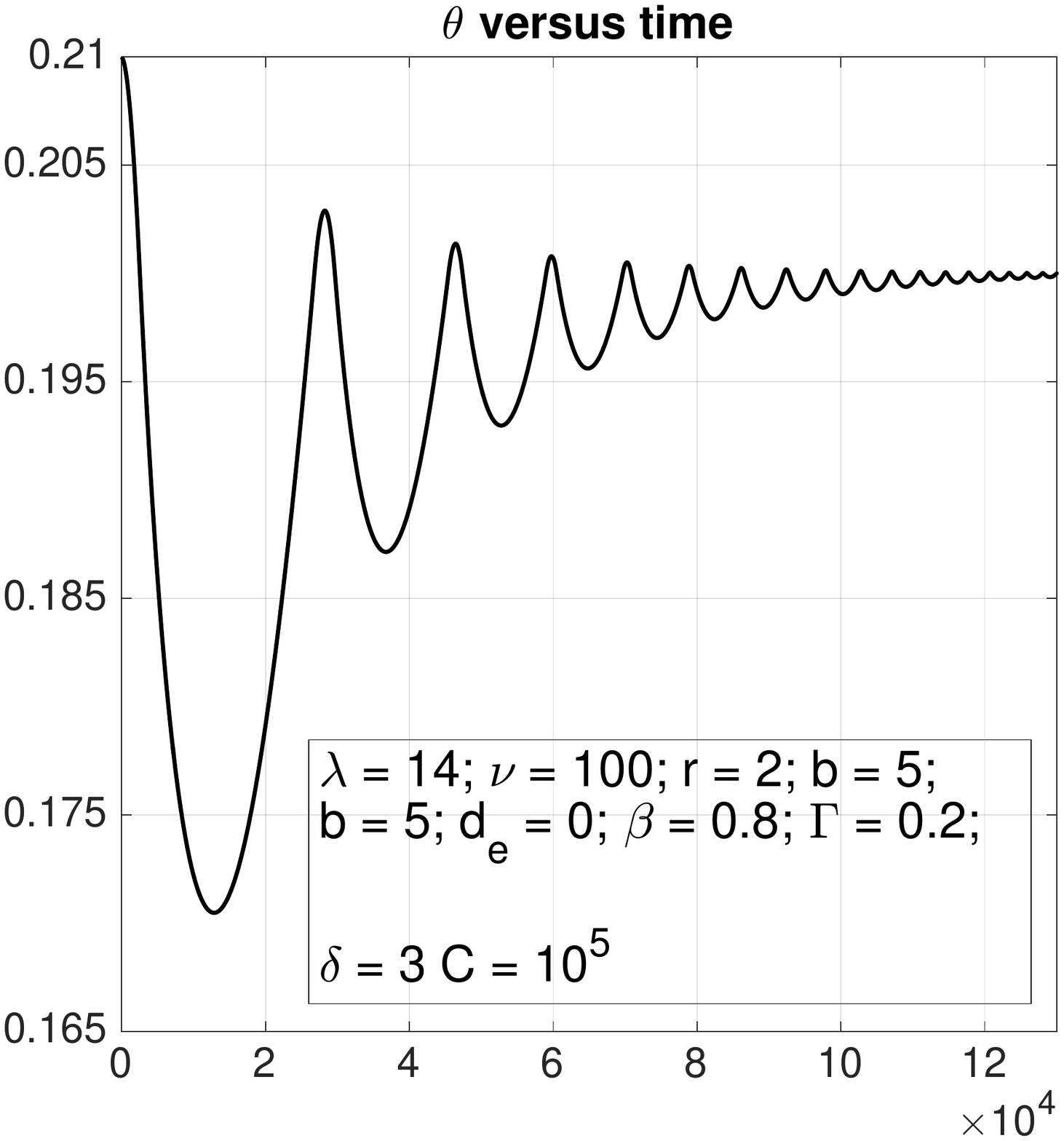}
    \vspace{-20mm}
    \caption{Theta trajectory}
    \label{fig_theta}  
\end{minipage}
\hspace{10mm}
\begin{minipage}{5cm}
      \centering
    \includegraphics[width=6cm, height = 8cm]{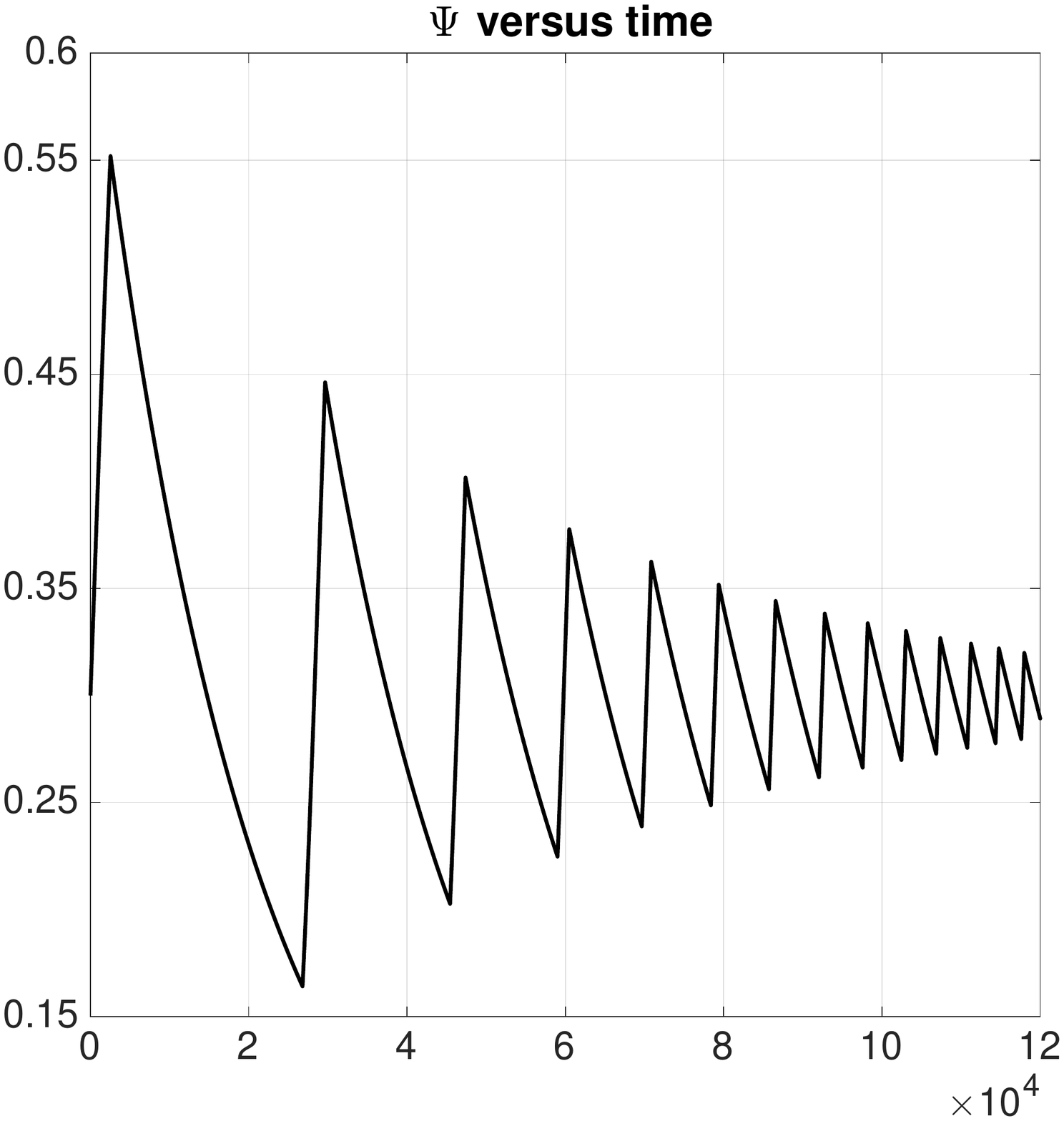}
      \vspace{-20mm}
    \caption{Caption}
    \label{fig_psi}  
\end{minipage}
\vspace{-10mm}
\end{figure}}

 \section{Conclusions}

With the ongoing pandemic in mind, we consider a scenario where  the vaccines are being prematurely introduced. Further, due to lack of information about the side-effects  and efficacy of  the vaccine, individuals exhibit {vaccination hesitancy}. This chaos is further amplified sometimes due to reported disease statistics, unavailability of the vaccine, leading to {vaccination urgency}. We developed an epidemic SIS model to capture such aspects, where the 
system changes due to births, deaths, infections and recoveries, while influenced by the dynamic vaccination decisions of the population. As observed in reality, a variety of behavioral patterns are considered, in particular follow-the-crowd, free-riding and vigilant agents.  Using stochastic approximation techniques, we derived the time asymptotic proportions of the infected and vaccinated population for a given vaccination response. Additionally,  we considered the conflict of stability for dynamic policies against static mutations, identified the strategies which are  stable against static mutations and studied the corresponding equilibrium states. 

Interestingly, the agents exhibit different behaviors and lead to different equilibrium states, however, at ESS-AS, all of the agents reach same limit state, where they choose vaccination either with probability $1$ or $0$, only based on system parameters. 
Also, by analysing the corresponding ODEs, we obtained many responses under which disease can be eradicated completely, but none of those are stable against mutations. Ironically, this is a resultant of the rationality exhibited by agents, which prevents them from reaching the disease-free state.

We observed certain surprising patterns at evolutionary stable equilibrium: (a) no one gets vaccinated with abundant vaccines; scarcity makes them rush; (b) the limit infected proportions   are concave functions of birth rates. 


\hide{Further, it is observed that the infected proportions achieved at ESS-AS are concave functions of birth rates, i.e., with high infection rate per birth, the population proactively vaccinates and infected fraction is lower, than that of lower infection rate per birth. With further increase in birthrate, the infected fraction decreases to zero.

Also, by analysing the corresponding ODEs of the system, we obtained many responses under which disease can be eradicated completely, but none of those are stable. Ironically, this is a resultant of rationality exhibited by agents, i.e., they become rational once the infection proportions reduces and reach the non-disease free state.  

}

Lastly,   the excess deaths did not change the patterns of ES equilibrium states versus parameters, however, the limit  vaccination fractions are much higher. So, in all it appears individuals rush for vaccine and we have smaller infected fractions at ES equilibrium, when there is a significant scare (of either deaths,  or of scarcity of vaccines, or high infection rates etc.). Ironically, the disease can be better curbed with excess deaths.

\section*{Appendix A: Stochastic approximation related proofs}
\begin{lemma}\label{lemma_error_convergence_2}
Let $\delta = 2/(N(0)-1)$. 
Then for any $k$,  $\eta_k \ge {\bar \delta} := \frac{N(0)-3}{(N(0)-1)^2}$ a.s.  And thus,  \vspace{-5mm}
$$E\left [ \frac{1}{ \eta_k^2} \right ] \le \frac{1}{{\bar \delta }^2} \mbox{ and }
E_k\left[\left|\left(\frac{1}{\eta_{k+1}}-\frac{1}{\eta_k}\right)\right|\right]\le \epsilon_k \frac{\bar{\delta}+1}{\bar{\delta}^2} \  \mbox{ a.s.,  for any } k.$$
\end{lemma}
\TR{Proof is provided in \cite{TR}. \eop}{\textbf{Proof:} 
We consider $\delta = 2/(N(0)-1)$, where $N(0)$ is the initial population. The minimum value that $\eta_{k+1}$ can take for any transition epoch  $k$ with  $\eta_k > \delta$ equals (this happens when death occurs, see \eqref{eqn_stoch_approx})
 $$
 \eta_k \left (1-\epsilon_k\right ) - \epsilon_k  \ge \delta (1-\epsilon_k) - \epsilon_k,
 $$
  and the first possible epoch $k$ at which $\eta_{k+1} $ drops below $\delta$,  is at $k = N(0)-2$  (this happens when the first   $N(0)-1$ transition epochs    are  all due to deaths),  and, hence the minimum possible value that $\eta_{k+1}$ can take with  $\eta_k > \delta$  equals
 \begin{eqnarray}
       {\bar \delta} := 
 \delta \left (1- \frac{ 1}{N(0)-1} \right ) - \frac{1}{N(0)-1} = 
 \frac{N(0)-3}{(N(0)-1)^2}  > 0. \label{Eqn_delta_bar}
 \end{eqnarray}

From equation \eqref{Eqn_VNSk} and \eqref{eqn_stoch_approx},  we have (because  $|\eta_k-G_{N,k+1}| \le \eta_k + 1$):
\begin{eqnarray*}
     \eta_{k+1}-  \eta_k &=& \indc{ \eta_k > \delta}\epsilon_k\left( G_{N,k+1}-\eta_k \right) \\
     \implies \frac{1}{\eta_{k+1}}-\frac{1}{\eta_k}   & = &   \frac{ \indc{ \eta_k >  \delta}   \epsilon_k }{\eta_{k+1}\eta_k}(\eta_k-G_{N,k+1}), \mbox{ and so, } 
 \\ 
     \left| \frac{1}{\eta_{k+1}}-\frac{1}{\eta_k}\right| &  \le &  \epsilon_k\frac{\bar{\delta}+1}{\bar{\delta}^2} \  \mbox{ a.s.} 
     \hspace{5.4cm}
     \mbox{ \eop }
\end{eqnarray*}}

\begin{lemma}\label{Lemma_alpha_km}
The term $\alpha^m_k  \to 0$ a.s.,   and,     
$\sum_k \epsilon_k |\alpha^m_k|<\infty$ a.s. for $m=\theta, \psi$.
\end{lemma}
\textbf{Proof:} We will provide the proof 
for $\alpha^{\psi}_k$ and proof goes through in exactly similar line  for $\alpha^\theta_k$.  From equation \eqref{eqn_Lk}, as in  \eqref{en_L_{k+1}},

\vspace{-6mm}
{\small 
$$ \alpha^{\psi}_k = E_k \left [ L_{k+1}^\psi - \frac{\eta_{k+1} }{\eta_{k}} L_{k+1}^\psi  \right ], \mbox{ where } L_{k+1}^\psi=\frac{\indc{\eta_k > \delta}}{\eta_{k+1}}\left[ G_{V,k+1} -   (N_{k+1}- N_k) \psi_k  \right].$$} 
By Lemma \ref{lemma_error_convergence_2} 
and because $| G_{V,k+1} -   (N_{k+1}- N_k) \psi_k | \le 2$ a.s., we have:

\vspace{-4mm}
{\small 
\begin{eqnarray}\label{eqn_lemma3}
 \left | \alpha^{\psi}_k  \right | \TR{\  \le \ }{ &\leq& E_k \left [ \left |  \left(\frac{1}{\eta_{k+1}}-\frac{1}{\eta_k}\right)\left[ G_{V,k+1} -   (N_{k+1}- N_k) \psi_k  \right] \right | \right ],\nonumber\\ 
 &\le & }
  \ 2E_k \left [  \left|\frac{1}{\eta_{k+1}}-\frac{1}{\eta_k}\right| \right ] \le   2 \epsilon_k\frac{\bar{\delta}+1}{\bar{\delta}^2} \mbox{ a.s.}
\end{eqnarray}}
Thus we have: \vspace{-8mm}
\begin{eqnarray*} \hspace{15mm}
\sum_{k=1}^{\infty}\epsilon_k |\alpha_k^\psi| &\le&2\frac{\bar{\delta}+1}{\bar{\delta}^2} \sum_{k=1}^{\infty}\epsilon^2_k  <\infty  \ \  a.s.  \hspace{3cm} \mbox{ \eop}
  \end{eqnarray*} 
  \vspace{-4mm}


 \begin{lemma}\label{lem_sum_lk}
 $\sup_k E|L^m_k|^2<\infty$ for $m=\theta, \psi,\eta$.
 \end{lemma}
 \textbf{Proof:}
 The result follows by Lemma \ref{lemma_error_convergence_2} and \eqref{eqn_Lk} (for an appropriate $C$):
\begin{eqnarray*}
|L^\theta_k|^2 \  \le \ \frac{4}{\eta^2_{k}} \ a.s., \mbox{ and}\ |L^\psi_k|^2 \  \le \ \frac{4}{\eta^2_{k}} \ a.s., \mbox{ and} \ |L^\eta_k|^2\le C<\infty, \ \ a.s. \hspace{1.4cm}
\mbox{ \eop}
\end{eqnarray*}


\noindent {\bf Proof of Theorem \ref{thrm1}:} 
As in \cite{kushner2003stochastic}, we will show that the 
following sequence  of piece-wise constant functions that  start with $\up_k$  are equicontinuous in extended sense.  Then the result follows from  \cite[Chapter 5, Theorem 2.2]{kushner2003stochastic}. 
  Define $(\up^k(t))_k:=( \theta^k(t),\psi^k(t), \eta^k(t))_k$ where,
 
 \vspace{-5mm}
 {\small \begin{eqnarray*}
          \theta^k(t)=\theta_k+\hspace{-2mm} \sum_{i=k}^{m(t_k+t)-1}\hspace{-5mm} \epsilon_k L^{\theta}_{k+1},  \  \ 
           \psi^k(t)=\psi_k+\hspace{-2mm} \sum_{i=k}^{m(t_k+t)-1}\hspace{-5mm} \epsilon_k L^{\psi}_{k+1},\  \ 
            \eta^k(t))_k=\eta_k+\hspace{-2mm} \sum_{i=k}^{m(t_k+t)-1}\hspace{-5mm} \epsilon_k L^{\eta}_{k+1}, 
\end{eqnarray*}}where $m(t):=\max\{k:t_k\le t\}$.
 This proof is exactly similar to that provided in the proof of \cite[Chapter 5, Theorem 2.1]{kushner2003stochastic} for the case with continuous $g$, except for the fact that $g (\cdot)$ in our case is not continuous. We will only provide differences in the proof steps towards 
       $(\theta^k(t))_k$ sequence, and it can be proved analogously for others.
      \TR{ Towards this we define $M_k^{\theta}=\sum_{i=0}^{k-1}\epsilon_i\delta M^{\theta}_i $ with
       $\delta M^{\theta}_k :=  L^{\theta}_{k+1}-g^{\theta}(\up_k)-\alpha^\theta_k$
       as in \cite{kushner2003stochastic}
       and show the required uniform continuity properties  in view of Lemmas \ref{Lemma_alpha_km} - \ref{lem_sum_lk}. Observe 
       that 
       $\eta_k \le 1+ N(0)/k$,  
        $\theta_k \le  1$ for any $k$. Now the uniform continuity of integral terms like the following is achieved because our $g (\cdot) $ are bounded:
       
       \vspace{-5mm}
       {\small 
       \begin{eqnarray*}
     \left|\int_s^t g^\theta(\up^k(z)) dz \right| \le \int_s^t|g^\theta(\up^k(z)) |dz\le\int_s^t\frac{\lambda+r+b+2 d_e}{\delta(d+b)}dz 
        \le  \bar{m}(t-s)\le \bar{m}\delta_1.
\end{eqnarray*}} 
 Such arguments lead to the required equicontinuity (details are in \cite{TR}). \eop 
 }{Towards this, define $\delta M^{\theta}_k :=  L^{\theta}_{k+1}-g^{\theta}(\up_k)-\alpha^\theta_k$ and then, equation \eqref{eqn_stc_theta} can be re-written as,
     \begin{eqnarray*}
          \theta_{k+1}=\theta_k+\epsilon_k \delta M^{\theta}_k+ \epsilon_k g^{\theta}(\up_k)+\epsilon_k \alpha_k^\theta,
     \end{eqnarray*} which gives us:
     \begin{eqnarray*}
     \theta^k(t)=\theta_k+ \sum_{i=k}^{m(t_k+t)-1} \epsilon_k\delta M^{\theta}_k  +\sum_{i=k}^{m(t_k+t)-1} \epsilon_k g^{\theta}(\up_k)+\sum_{i=k}^{m(t_k+t)-1} \epsilon_k\alpha_k^\theta.
       \end{eqnarray*}Now define $M_k^{\theta}=\sum_{i=0}^{k-1}\epsilon_i\delta M^{\theta}_i $, it is easy to prove $(M_k^{\theta},\mathcal{F}_k)_k$ is martingale, where $\mathcal{F}_k$ is natural filtration. Thus, using Martingale inequality (see \cite[ chapter 4, equation (1.4)]{kushner2003stochastic} for $q(M)=M^2$), we get for each $\mu>0$
       \begin{equation*}
           P\left\{\sup_{i\le j \le k}|M^{\theta}_j-M^{\theta}_i|\ge \mu\right\}\le \frac{E|\sum_{m=i}^{k-1}\epsilon_m\delta M^{\theta}_m|^2}{\mu^2}
       \end{equation*} Using the fact $E[\delta M^{\theta}_i \delta M^{\theta}_j
       ]=0$ for $i<j$,       
       \hide{Using this we get,
       
       \vspace{-2mm}
       {\small
       \begin{eqnarray*}
           P\left\{\sup_{i\le j \le k}|M^{\theta}_j-M^{\theta}_i|\ge \mu\right\}&\le& \frac{\sum_{m=i}^{k-1}\epsilon^2_mE|\delta M^{\theta}_m|^2}{\mu^2} 
           \\
           &\le& \frac{2
           \sum_{m=i}^{k-1}\epsilon^2_m \left(E[ L^{\theta}_{m+1}-g^\theta(\up_m)]^2+(\alpha^\theta_m)^2\right)}{\mu^2}
            \end{eqnarray*}}} Lemma \ref{lem_sum_lk} and equations  \eqref{eqn_ODE} and \eqref{eqn_lemma3}, we have: \hide{for an appropriate constant $K < \infty$:  $\sup_k E[ L^{\theta}_{m+1}-g^\theta(\up_m)]^2+(\alpha^\theta_m)^2=K<\infty$ for some constant $K$. Thus we have:}
           \begin{eqnarray}\label{eqn_martingale}
           \hide{P\left\{\sup_{i\le j \le k}|M^{\theta}_j-M^{\theta}_i|\ge \mu\right\}  &\le& \frac{K
           \sum_{m=i}^{k-1}\epsilon^2_m}{\mu^2} \le \frac{K}{\mu^2}         \sum_{m=i}^{\infty}\epsilon^2_m \mbox{ and so }\nonumber\\}
  \lim_{i\rightarrow\infty}  P\left\{\sup_{i\le j }|M^{\theta}_j-M^{\theta}_i| \geq \mu \right\}&=&0 \mbox{ for each } \mu>0.
               \end{eqnarray}   
        Now, we can re-write $\theta^k(t)$ as
        \begin{equation*}
            \theta^k(t)=\theta_k+\int_0^t g^\theta(\up^k(z)) dz + M^{k,\theta} (t) + A^{k,\theta} (t) + \rho^{k,\theta} (t) 
        \end{equation*}
        where we denote $ M^{k,\theta} (t)= \sum_{i=k}^{m(t_k+t)-1} \epsilon_i\delta M^{\theta}_i $, and  $A^{k,\theta} (t) =\sum_{i=k}^{m(t_k+t)-1} \epsilon_i\alpha_i^\theta $ and $ \rho^{k,\theta} (t) =\sum_{i=k}^{m(t_k+t)-1} \epsilon_i g^\theta(\up_i)-\int_0^t g^\theta(\up^i(z)) dz$. 
        
        Clearly $\theta^k(0)=\theta_k \le 1$ and thus to claim equi-continuity, we need: for each $T$ and $\epsilon>0$, there is a $\delta>0$ such that
        \begin{eqnarray*}
             \lim\sup_n\sup_{0\le t-s\le \delta,|t|<T }|\theta^k(t)-\theta^k(s)|<\epsilon.
        \end{eqnarray*}
        Towards this, we have,
        \hide{
\begin{eqnarray*}
   |\theta^k(t)-\theta^k(s)|   &\le&\sup\left|\int_s^t g^\theta(\up^k(z)) dz\right|  +\sup\left| M^{k,\theta} (t) - M^{k,\theta} (s) \right|\\
   &&\hspace{2cm}+\sup\left|A^{k,\theta} (t)-A^{k,\theta} (s)\right|+\sup\left| \rho^{k,\theta} (t)- \rho^{k,\theta} (s)\right|
   \end{eqnarray*}and so, }
   \begin{eqnarray*}
     \sup   |\theta^k(t)-\theta^k(s)|&\le&\sup\left|\int_s^tg_k(\theta^n(z)) dz\right|  +\sup\left| M^{k,\theta} (t) - M^{k,\theta} (s) \right|\\   
   &&\hspace{2cm}+\sup\left|A^{k,\theta} (t)-A^{k,\theta} (s)\right|+\sup\left| \rho^{k,\theta} (t)- \rho^{k,\theta} (s)\right|
   \end{eqnarray*}
where supremum is taken over $S_T := \{ 0\le t-s\le \delta,|t|<T\}$. Let us consider the first term from above, using equation \eqref{eqn_ODE}, with $\bar{m} := \frac{\lambda+r+b+2 d_e}{\delta(b+d)}$:
\begin{eqnarray*}
     \left|\int_s^t g^\theta(\up^k(z)) dz \right| &\le& \int_s^t|g^\theta(\up^k(z)) |dz 
       \le \bar{m}(t-s)\le \bar{m}\delta_1, \mbox{ for some } \delta_1 > 0. 
\end{eqnarray*}
For the second term, from \eqref{eqn_martingale}, by continuity of probability,
\begin{eqnarray*}
        P\left\{\lim_{i\rightarrow\infty}\sup_{i\le j }|M^{\theta}_j-M^{\theta}_i|\ge\mu\right\}&=&0
\end{eqnarray*}
Let $A_n:=\left\{\omega: \lim_{i\rightarrow\infty}\sup_{i\le j }|M^{\theta}_j-M^{\theta}_i|<\frac{1}{n}\right\}$, then $P(A_n)=1$ for each $n>0$. Now we  aim to show that $(M^{k,\theta}(\omega,\cdot)) $ converges to zero uniformly on each bounded interval in $(0,\infty)$ as $k\rightarrow\infty$, for each $\omega \not\in N:=\Omega-\cap_{n>0}A_n$, i.e., $\sup_{S_T} M^{k,\theta}(t)\rightarrow 0$. To this end, we consider $S_T$, then for every $\omega\not\in N$,
\hide{
\begin{eqnarray*}
     |M^{k,\theta}(t)|&=&\left|\sum_{i=k}^{m(t_k+t)-1} \epsilon_i\delta M^{\theta}_i\right|=\left|\sum_{i=0}^{m(t_k+t)-1} \epsilon_i\delta M^{\theta}_i-\sum_{i=0}^{k-1} \epsilon_i\delta M^{\theta}_i\right|=|M^{\theta}_{m(t_k+t)}-M^{\theta}_k|
\end{eqnarray*}This gives us }
\begin{eqnarray*}
     \sup_{t \in S_T}|M^{k,\theta}(t)|&=& \sup_{t \in S_T }|M^{\theta}_{m(t_k+t)}-M^{\theta}_k|\le\sup_{m(t_k+t)\ge k}|M^{\theta}_{m(t_k+t)}-M^{\theta}_k|.
\end{eqnarray*}\hide{and so
\begin{eqnarray*}
\lim_{k\rightarrow\infty} \sup_{t\in S_T}|M^{k,\theta}(t)|&\le&\lim_{k\rightarrow\infty}\sup_{m(t_k+t)\ge k}|M^{\theta}_{m(t_k+t)}-M^{\theta}_k|<\frac{1}{n}
\end{eqnarray*}
taking $n\rightarrow \infty$,} Then, taking $k \to \infty$ first and then letting $n$ go to $\infty$, we get our claim (see definition of $A_n$). For the third term, uniform convergence to zero follows from Lemma \ref{Lemma_alpha_km} as shown below:
\begin{eqnarray*}
\sup_{S_T}|A^{k,\theta} (t)| &\leq&\sup_{S_T}\sum_{i=k}^{m(t_k+t)-1}|\epsilon_i\alpha_i^\theta|\le \sum_{i=k}^{\infty}\epsilon_i|\alpha_i^\theta|,
\end{eqnarray*} and taking $k\rightarrow \infty$, we get the result.
\hide{For the last term, first we will prove for $t=t_n-t_k, \ (n>k)$ that $\rho^{k,\theta}(t)=0$. Towards this, let $n=k+1$, we have $t_n-t_k=\epsilon_k$:
\begin{eqnarray*}
  \rho^{k,\theta} (t_n-t_k) =\sum_{i=k}^{m(t_{k+1})-1} \epsilon_k g^\theta(\up_k)-\int_0^{\epsilon_k} g^\theta(\up^k(z)) dz=0.
\end{eqnarray*}
because, for $0 \le s\le \epsilon_k$, $\theta(\up^k(z))=\up_k$. Let the claim be true for $n=k+2,\dots,k+l$, then, for $n=k+l+1$, we have:
\begin{eqnarray*}
  \rho^{k,\theta} (t_{k+l+1}-t_k) =\rho^{k,\theta} (t_{k+l}-t_k)+ \epsilon_{k+l} g^\theta(\up_{k+l})-\int_0^{\epsilon_{k+l}} g^\theta(\up^k(z)) dz=0.
\end{eqnarray*}
Thus, by induction the claim is true.} For the last term, it can be proved by induction that when $t$ exactly corresponds to the end of epochs, i.e., when $t = t_n-t_k, \ (n>k)$ that $\rho^{k,\theta}(t)=0$.  Now, we are only left to prove that $\rho^{k,\theta}(t)\rightarrow 0$ uniformly in $t$ (for general $t$) as $k \rightarrow \infty$\hide{, i.e., $|\rho^{k,\theta}(t)-0|<\epsilon$ for every $k>N_{\epsilon}$}. We will prove this claim for each $T>0$, such that $|t|<T$:

\vspace{-2mm}
{\small
\begin{eqnarray*}
 |\rho^{k,\theta}(t)| &=&\left| \sum_{i=k}^{m(t_k+t)-1} \hspace{-4mm} \epsilon_i g^\theta(\up_i)-\int_0^t g^\theta(\up^k(z)) dz\right| = \left|\int_{m(t_k+t)}^t   g^\theta(\up^k(z)) dz\right| \\
 &\le& (t - m(t_k + t))\bar{m} \leq \epsilon_k \bar{m},
\end{eqnarray*}}where $\epsilon_k \to 0$ and first inequality follows as in proof of first term. This proves the equicontinuity in extended sense for $(\theta^k(z))_k$. Proof follows in exact similar lines for $(\psi^k(z))_k$. For $(\eta^k(z))_k$ we need to prove  $\eta^k(0)=\eta_k$ is bounded. From \eqref{Eqn_VNSk},
\begin{eqnarray*}
 \eta_{k}&\le& \frac{1}{k}\left(N(0)+\sum_{i=0}^{k} \mathbb{B}_{i}\right) \le \frac{1}{k}\left(N(0)+k\right)= 1+\frac{N(0)}{k}.
\end{eqnarray*} Rest of the proof follows in the same way as above.  This proves $(\up^k(t))_k$ is equicontinuous in extended sense. From \cite[Chapter 5, Theorem 2.2]{kushner2003stochastic} $\up_k\rightarrow A$.

\noindent\textbf{Part (ii):} From equicontinuity of $(\up^k(\cdot))$ and extended version of Arzela-Ascoli Theorem \cite[Chapter 4, Theorem 2.2]{kushner2003stochastic},for $\omega\not\in N$, there is a sub-sequence $(\up^{k_m})$ which converges to some continuous  limit uniformly on each bounded interval. It is easy to verify that the limit satisfies the ODE \eqref{eqn_ODE}. Basically,  $(\up^{k_m}(\cdot))$  converges to the solution of \eqref{eqn_ODE} uniformly on each bounded interval. Denote this solution as $\up_{*}(\cdot)$  Now, $(\up^{k_m}(0))=(\up_{k_m})$ converges to $\up_{*}(0)$, which implies  $\up_{*}(0)=\lim_{k_m}\up_{k_m}$. Further from uniform convergence, for any $\delta>0$, there exists a $\bar{N}$ such that for any $T\in(0,\infty)$,
$$\sup_{0\le t\le T}d(\up^{k_m}(t),\up(t))<\delta\mbox{ for all } k_m>\bar{N}.$$ Consider $t=t_k-t_{k_m}, (k>k_m)$ such that $0\le t\le T$, i.e., $t_{k_m}\le t_k\le t_{k_m}+T$. Note that $\up^{k_m}(t_k-t_{k_m})=\up_k$, so we have for all $k_m>\bar{N}$
 $$
            \sup_{k: t_k \in [t_{k_m}, t_{k_m} + T]} d(\up_k, \up_{*}(t_k - t_{k_m}))<\delta,
            $$
which proves the claim. Observe here that the sub-sequence only depends upon $\omega$ and not on $T$. \eop}
\section*{Appendix B: ODE attractors related proofs}

%
\newcommand{\ttheta}{\widetilde{\theta}}
\newcommand{\tpsi}{\widetilde{\psi}}
%

\noindent \textbf{Proof of Theorem \ref{thrm_FCagents}:} 
Let $\Upsilon := (\theta, \psi, \eta)$. Let $\hat{\Upsilon}$ represent the  corresponding attractors from Table \ref{table_FC}. Here  $q(\theta, \psi) = \min\{\tilde{q}(\theta, \psi), 1\} $ with $\tilde{q}(\theta, \psi) = \beta \psi$.
    
    We first consider the case where $\tilde{q}(\htheta, \hpsi) < 1$.
Further, note that one can re-write ODEs,
$\dot{\up} = g(\up)$, as below:
    $$
    \dot{\theta} = \frac{\indc{\eta >\delta}A \theta}{\eta \varrho }, \  \dot{\psi} = \frac{\indc{\eta >\delta}B  \psi}{\eta \varrho }, \mbox{ and } \dot{\eta} =  \indc{\eta>\delta} C,  \mbox{ where }
    $$
    $A = A(\Upsilon) := (1- \theta - \psi)\lambda - r - b$,  $B = B(\Upsilon) := (1- \theta - \psi)\beta \nu - b$ and $C = C(\Upsilon)  := \nicefrac{(b-d)}{\varrho} - \eta$.  To this end, we define the following Lyapunov function based on the regimes of parameters:
    
    \vspace{-2mm}
    {\small 
\[
V(\Upsilon) := 
\begin{cases}
\left(\hat{A}(\theta)\right)^2 +  \psi \left(\hat{B}(\psi)\right)^2 + C(\eta) (\hat{\eta} - \eta),  &\mbox{ if } \beta < \rho \mu, \rho > 1, \\
\theta \left(\hat{A}(\theta)\right)^2 +   \left(\hat{B}(\psi)\right)^2 + C (\eta) (\hat{\eta} - \eta), &\mbox{ if } \beta > \rho \mu,  \mbox{ and }  \rho > 1\\
\theta \left(\hat{A}(\theta)\right)^2 +  \psi \left(\hat{B}(\psi)\right)^2 + C (\eta)(\hat{\eta} - \eta), &\mbox{ if } \rho < 1,
\end{cases}
\]}where   $\hat{A}(\theta) := A(\theta, \hat{\psi}, \hat{\eta})$,  $\hat{B}(\psi) := B(\htheta, \psi, \hat{\eta})$. 
\TR{ We complete this proof using the above functions and the details are in \cite{TR}. \eop

}{
Observe that $V$ is continuously differentiable, 
$V(\hat{\Upsilon}) = 0$ and   $V(\Upsilon) > 0$ for all $\Upsilon \ne \hat{\Upsilon}$. In all the computations below we consider a neighbourhood of $\hat{\Upsilon}$ for which $\eta > \delta$ and omit this indicator in further computations.
    
We begin with the case when $\beta < \rho \mu, \rho > 1$  Then, the derivative of $V(\up(t))$ with respect to time is given by:
    \begin{align}\label{eqn_derv_V_FC}
     \dot{V} &= \left<\nabla V, g(\up) \right> \\
     &= \left(-2\lambda \hat{A}(\theta) A \theta + \Big((\hat{B}(\psi))^2 - 2\beta \nu \psi  \hat{B}(\psi) \Big) B \psi  \right) \frac{1}{\eta \varrho }- (C(\eta) + \hat{\eta} - \eta)C. \nonumber
    \end{align}
    Since $\hat{B}(\hpsi) < 0$, by continuity of $B(\cdot)$, there exists a neighborhood of $\hat{\up}$ such that  in the neighborhood (with $\Delta > 0$) following is true:
    
    \vspace{-6mm}
    \begin{align}\label{eqn_bound_FC}
        B(\up) < 0, \hat{B}(\psi) < 0 \mbox{ and } \frac{1}{\eta \varrho }(\hat{B}(\psi))^2   B(\up) < -\Delta.
    \end{align}
     Observe  $\hat{A}(\theta)A = (\hat{A}(\theta))^2 - \psi \lambda^2 (\hat{\theta} - \theta)$. Using continuity arguments again,
     choose a neighborhood of $\hat{\up}$ (further smaller, if required) such that  (see \eqref{eqn_bound_FC})
    $$
    \frac{\psi}{\eta \varrho}\left[2 \theta  \lambda^3 (\htheta - \theta) +  (\hat{B}(\psi))^2   B(\theta, \psi)   \right]  < 0.
    $$
    This proves that the first component in \eqref{eqn_derv_V_FC} is negative. We are now left to prove that second component, $- (C + \hat{\eta}- \eta) C$, is also negative. To this end, notice $(\heta- \eta) C = (\heta - \eta) (\nicefrac{(b-d)}{\varrho} - \eta) > 0$  (observe $\heta = \nicefrac{b-d}{\hat{\varrho}}$, and $(\heta- \eta) C(\heta) > 0$), which proves the claim for second component as well. 
    Conclusively, we get that $\dot{V} < 0$ in the neighborhood chosen above.
    This implies that $\hat{\up}$ is locally asymptotically stable in the sense of Lyapunov, i.e., $\up \to \hat{\up}$ when the ODEs start in the neighborhood of $\hat{\up}$.  Further, the arguments follow in exactly similar lines in rest of the cases. Lastly,   when $\tilde{q}(\htheta, \hpsi) > 1$ and $(\htheta, \hpsi)$ is a boundary point, then the ODE corresponding to $\psi$ changes to :
    $$
     \dot{\psi} = \frac{\indc{\eta >\delta}B }{\eta \varrho },  \mbox{ where }B = B(\Upsilon) := (1- \theta - \psi) \nu  - b \psi.
    $$
    Then, the proof follows exactly as above (note here $\hat{A}(\htheta) < 0$, since $\mu \rho < \mu + 1$). For the interior point, proof follows as in Lemma \ref{lemma_ODE}.  \eop} 
  
  \TR{}{  
\vspace{2mm}
\noindent \textbf{
Proof of Theorems \ref{thrm_FR_agents}-\ref{thrm_VFC1}:} Let $\hat{\up}$  be the  corresponding attractors from Table \ref{table_FR}. When $\hat{\up}$ are boundary points,  the proof follows as shown in proof of Theorem \ref{thrm_FCagents}. For the interior point, $\hat{\up}$, an appropriate Lyapunov function can be constructed as described in Lemma \ref{lemma_ODE} to complete the proof.  \eop}

\begin{lemma}\label{lemma_ODE}
Let $\htheta, \hpsi > 0$. If $\tilde{q}(\htheta, \hpsi) \ne  1$,   there exists a Lyapunov function such that
$(\htheta, \hpsi, \hat{\eta})$  is  locally asymptotically stable attractor for ODE \eqref{eqn_ODE} in the sense of Lyapunov.  
\end{lemma}

{\bf Proof:} 
\TR{We use similar notations as in previous proof.}{
Let $\hat{\up}$ be the given interior attractor for the given $\pi \in \Pi$. One can re-write ODEs, $\dot{\up} = g(\up)$ as:
$$
\dot{\theta} = \frac{\lambda \indc{\eta >\delta}A \theta}{\eta \varrho }, \  \dot{\psi} = \frac{\nu \indc{\eta >\delta}B  }{\eta \varrho }, \ \dot{\eta} = \indc{\eta >\delta} C,  \mbox{ where }
$$
    $A = A(\Upsilon) := 1- \theta - \psi - \frac{1}{\rho}$,  $B = B(\Upsilon) := (1- \theta - \psi)q(\theta, \psi)  - \mu \psi$,  and $C = C(\Upsilon)  := \nicefrac{(b-d)}{\varrho} - \eta$.
In all the computations below, we consider a neighbourhood of $\hat{\Upsilon}$ for which $\eta > \delta$ and omit this indicator in further computations.

}
Let us first consider the case where $\tilde{q}(\htheta, \hpsi) < 1$, i.e., $q(\htheta, \hpsi) = \tilde{q}(\htheta, \hpsi)$. Then, one can choose a neighborhood (further smaller, if required) such that $\hat{q} - \delta < q(\theta, \psi) < \hat{q} + \delta$, and $q(\theta, \psi) = \tilde{q}(\theta, \psi)$ for some $\delta > 0$. Define the following Lyapunov function (for some $w_1, w_2 > 0$, which would be chosen appropriately later):
\begin{align}
    V(\up) := w_1 (\htheta - \theta) \hat{A}(\theta) + w_2 (\hpsi - \psi) \hat{B}(\psi) + C(\hat{\eta} - \eta), \mbox{ where }
\end{align}
$\hat{A}(\theta) := 1-\theta-\hpsi - \frac{1}{\rho}$, and $\hat{B}(\psi) := \hat{q}(1-\htheta-\psi) - \mu \psi$ (recall $\hat{q} := q(\htheta,\hpsi)$). 
\TR{Call }{Observe   $V$ is continuously differentiable, 
$V(\hat{\Upsilon}) = 0$ and   $V(\Upsilon) > 0$ for all $\Upsilon \ne \hat{\Upsilon}$.

For simplicity in notations, call} $\ttheta := \htheta - \theta$ and $\tpsi := \hpsi - \psi$.    
The derivative of $V(\up(t))$ with respect to time is:

\vspace{-2mm}
{\small 
    \begin{align}\label{eqn_derv_V_lemma}
     \dot{V} &= \left<\nabla V, g(\up) \right>  \\
     &= -\left(\hat{A}(\theta) + \ttheta \right)  \frac{A\theta \lambda w_1}{\eta \varrho}  - \left(\hat{B}(\theta) + \tpsi(\hat{q} + \mu) \right)  \frac{B\nu  w_2}{\eta \varrho} - (C(\eta) + \heta - \eta)C . \nonumber
    \end{align}}
    
One can prove that the  last component, i.e., $- (C + \heta- \eta) C$ is strictly negative in an appropriate neighborhood of $\hat{\up}$ as in proof of Theorem \ref{thrm_FCagents}.  Now, we proceed to prove that other terms in $\dot{V}$ (see \eqref{eqn_derv_V_lemma}) are also strictly negative in a neighborhood of $\hat{\up}$.

Consider the term\footnote{Observe that $\hat{A}(\theta) A = (\hat{A}(\htheta) + \ttheta ) A$, and $\hat{A}(\htheta) = 0$.} $\left(\hat{A}(\theta) + \ttheta \right)A$, call it $A_1$:

\vspace{-2mm}
{\small
\begin{align}\label{A_term}
    \begin{aligned}
A_1 & = 2\ttheta A = 2\left( \ttheta^2 + \ttheta \tpsi \right) = 2\left( \left( \ttheta c_1 + \frac{1}{2c_1}\tpsi\right)^2 + (1-c_1^2)\ttheta^2 - \frac{1}{4c_1^2} \tpsi^2  \right),
\end{aligned}
\end{align}}where $c_1$ will be chosen appropriately in later part of proof. \TR{ Similarly the term corresponding to $B$ is (details in \cite{TR}), \vspace{-4mm}}{Further, we have\footnote{{Note that $B = (1- \theta - \psi)(q(\theta, \psi) - \hat{q}) + (\ttheta + \tpsi) \hat{q} + \mu \tpsi + B(\htheta, \hpsi)$, and $B(\htheta, \hpsi) = 0$.}}:
\begin{align}\label{B_term}
    \begin{aligned}
    B_1 := \left(\hat{B}(\theta) + \tpsi(\hat{q} + \mu) \right)B &= 2\tpsi(\hat{q} + \mu)B\\
    &\hspace{-3cm}=2\tpsi(\hat{q} + \mu)(1-\theta-\psi)\big(q(\theta, \psi) - \hat{q} \big)  + 2\mu \tpsi^2(\hat{q} + \mu)\\
    &\hspace{-2cm}+ 2(\hat{q} + \mu)\hat{q}(\ttheta \tpsi + \tpsi^2).
\end{aligned}
    \end{align}
Note that the first term in \eqref{B_term} (call it $B_2$) can be written as, for some functions $p_1(\up)$ and $p_2(\up)$:
$$
B_2 = 2(\hat{q} + \mu)\Big[ p_1(\up) \tpsi^2 - p_2(\up)\tpsi \ttheta \Big], 
$$
such that for an appropriate neighborhood, we have: $p_1 - \bar{\delta} \leq p_1(\up) \leq p_1 + \bar{\delta}$ and $p_2 - \bar{\delta} \leq p_2(\up) \leq p_2 + \bar{\delta}$.
Now using $B_2$, one can re-write $B_1$  as:}

\vspace{-2mm}
{\small 
\begin{align}\label{eqn_B1_term}
    \begin{aligned}
    B_1 \TR{}{&=2(\hat{q} + \mu) \left [ \left ( \mu  +\hat{q} +  p_1(\up) \right ) \tpsi^2 + \hat{q}
    \left ( 1 - \frac{p_2(\up)}{\hat{q}} \right ) \tpsi\ttheta
    \right  ]
    \\}
    &= 2 (\hat{q} + \mu) \left[ \tpsi^2\Big( \mu + p_1(\up) + \hat{q} (1-c_2^2)\Big)  -  \frac{1}{4c_2^2 \hat{q}} \left(\hat{q} - p_2(\up) \right)^2 \ttheta^2 \right]\\
&\hspace{1cm}+ 2(\hat{q} + \mu) \hat{q}  \left(c_2 \tpsi + \frac{1}{2c_2} \left(1 - \frac{p_2(\up)}{\hat{q}} \right) \ttheta \right)^2.
\end{aligned}
\end{align}}
Thus, we get: $\dot{V} 
     <  -A_1 \frac{  \theta \lambda w_1}{\eta \varrho}  - B_1 \frac{\nu w_2}{\eta \varrho} $\TR{.}{(recall the last component in \eqref{eqn_derv_V_lemma} is strictly negative).}
Now, for $\dot{V}$ to be negative, we need (using terms, in \eqref{A_term} and \eqref{eqn_B1_term}, corresponding to $\ttheta^2, \tpsi^2$):

\vspace{-2mm}
{\small 
\begin{align*}
     \nu (\hat{q} + \mu)  &\leq \frac{4 c_2^2  \hat{q} }{ w_2 \left( \hat{q}- p_2(\up) \right)^2 }  \lambda w_1 (1-c_1^2) \theta, \mbox{ and }\\
     \frac{1}{2c_1^2}\theta \lambda w_1 &\leq  2\nu w_2 (\hat{q} + \mu)  \Big( \mu  + \hat{q} (1-c_2^2) + p_1(\up) \Big).
\end{align*}}
\TR{By appropriately choosing the constants (for various agents), we complete the proof (details in \cite{TR}). \eop}{
To this end, it is sufficient to choose a further smaller neighborhood of $\hat{\up}$ (call it again $\delta$-neighborhood, and such that $\htheta - \delta >0$) and some of the parameters $c_1, c_2, w_1$ and $w_2$  such that $1 - c_1^2, 1 - c_2^2 > 0$ and following is true: 
\begin{align}\label{eqn_bound}
    \nu (\hat{q} + \mu)   &= \frac{4 c_2^2\hat{q}}{ w_2 \left( \hat{q}- (p_2 - \delta) \right)^2 }  \lambda w_1 (1-c_1^2) (\htheta - \delta), \mbox{ and } \nonumber \\
    \frac{1}{2c_1^2}(\htheta + \delta) \lambda w_1 &\leq  2\nu w_2 (\hat{q} + \mu)  \Big( \mu  + \hat{q} (1-c_2^2) + p_1 - \delta \Big).
\end{align}
In above, first equation holds when $\hat{q} > p_2(\up)$, otherwise, the denominator changes to $w_2 \left( \hat{q}- (p_2 + \delta) \right)^2$.

Let $c_1^2 := \Omega_1 w_1$, and $c_2^2 := \Omega_2 w_2$, for appropriate $\Omega_1, \Omega_2$.
Then, by substituting the first equality in the second inequality of  \eqref{eqn_bound}: 

\vspace{-2mm}
{\small
\begin{align*}
    (\htheta + \delta)  &\leq \frac{ 2c_1^2    4 c_2^2\hat{q}}{ w_2 \left( \hat{q}- (p_2 - \delta) \right)^2 }   (1-c_1^2) (\htheta - \delta) 
    2 w_2    \Big( \mu  + \hat{q} (1-c_2^2) + p_1 - \delta \Big)
    \\
    &=     (\htheta - \delta) \Omega_1 w_1 (1- \Omega_1 w_1) \Omega_2 w_2 \Big( 1  + \frac{p_1 - \delta + \mu}{\hat{q}} -\Omega_2 w_2 \Big)\left( \frac{  4\hat{q} }{\hat{q} - (p_2-\delta)} \right)^2 \\
    &\leq    (\htheta - \delta)  \frac{1}{4}  
    \frac{1}{4}\Big( 1+ \frac{p_1 - \delta + \mu}{\hat{q}}   \Big)^2 \left( \frac{  4\hat{q} }{\hat{q} - (p_2-\delta)} \right)^2  \ =     (\htheta - \delta)   \left ( \frac{  1+ \frac{p_1 - \delta + \mu}{\hat{q}}  }{1 - \frac{ p_2-\delta}{\hat{q}}} \right )^2 , 
\end{align*}}
The last step is obtained by choosing   $\Omega_1w_1=1/2$ that provides the maximum value $1/4$ (for product  $\Omega_1w_1 (1-\Omega_1w_1)$) and $\Omega_2 w_2$ that provides the maximum value $(1+ (p_1-\delta+\mu)/ \hat{q} )^2/4 $ for term $\Omega_2 w_2 \Big( 1  + \frac{p_1 - \delta + \mu}{\hat{q}} -\Omega_2 w_2 \Big)$. 

Below we show that the above inequality can be solved for FC, FR as well as VFC agents. Basically we need to show the following 
\begin{align}\label{eqn_cond_lemma}
\frac{ \hat{q}+ p_1 + \mu }{ \hat{q} - p_2}  > 1, \mbox{ or, }  < -1.
\end{align}
To begin with, when $\hat{q}=1$ with $\tilde{q} (\htheta, \hpsi) > 1$, then one can chose a neighbourhood of the zero at which 
$\tilde{q} (\theta, \psi) > 1$, which implies $q = 1$ in that neighbourhood. Thus the above condition is true with $p_1 = 0, p_2 = 0$. This is true for all $\pi \in \Pi$. We will now prove that \eqref{eqn_cond_lemma} is true when $\tilde{q} (\htheta, \hpsi) < 1$.

For $\pi = FR$, using $2\tpsi(\hat{q} + \mu)(1-\theta-\psi)\big(q(\theta, \psi) - \hat{q} \big) $ in \eqref{B_term} and $B_2$, we have $p_1 = p_1(\hat{\up}) = -\hbeta (1-\htheta - \hpsi) (1-\hpsi - \psi)$ and $p_2 = 0$. For these $p_1, p_2$, equation \eqref{eqn_cond_lemma} is satisfied (since $\hbeta > \mu \rho$).

Proceeding as in $\pi=FR$ above, for $\pi = VFC1$ we get $p_1 = -(\hbeta \htheta)/\rho$ and $p_2 = (\hbeta \hpsi)/\rho$. Now, we have:
\begin{align*}
\frac{ \hat{q}+ p_1 + \mu }{ \hat{q} - p_2}  = \frac{\mu \rho^2 }{\mu \rho^2 - \hbeta} = 1 + \frac{\hbeta}{\mu \rho^2 - \hbeta}.
\end{align*}
If $\mu \rho^2 > \hbeta$, then, \eqref{eqn_cond_lemma} is true. It is also true when $\mu \rho^2 < \hbeta$ and then $2\mu \rho^2 > \hbeta$, as then the fraction is less than $-1$. In all the result is true for all $\hbeta \le 2 \mu \rho^2  $.
\eop}

\vspace{-4mm}
\section*{Appendix C: ESS related proofs}
\vspace{-1mm}
\TR{}{
\noindent\textbf{Proof of Lemma \ref{lem_ess_zero}:}
From Theorems \ref{thrm_FCagents}-\ref{thrm_VFC1} with $\rho>1$,  the attractors $(\hat{\theta},\hat{\psi})$ are $(1 - \frac{1}{\rho},0)$ for all  $\hat{\beta} \le {\bar \beta}_\pi$, with ${\bar \beta}_\pi> 0$.  For example, for $\pi = V_{FC}^1$, $\bar{\beta}_\pi =\nicefrac{\mu \rho^2}{\rho-1}$. 
%
%
If $h_m>0$,   from \eqref{eqn_utility}  the static-best response set   equals
$
{\cal B} (\pi(\hat {\beta} ) ) = \{0\} ,
$
and observe $q^*(0)=0$. 
Further by Lemma \ref{lem_cont_of_attract} given in this appendix (below),  best response against
$\epsilon$ mutational strategy, $\pi_\epsilon(\hat{\beta}, p)$ (with $\rho >1$) also equals $\{ 0\}$ for all $\epsilon \le {\bar \epsilon} $. 
Thus $\pi(\hat{\beta})$ with $\hat{\beta} = 0$ is an ESS-AS for any $\pi \in \Pi$, when $\rho >1$. 
 
From Theorems \ref{thrm_FCagents}-\ref{thrm_VFC1} with  $\rho<1$,  the attractors $(\hat{\theta},\hat{\psi})$ are $(0,0)$ for all $\hbeta$.   From equation \eqref{eqn_utility}, clearly $h(0,0) = c_{v_1} +{\bar c}_{v_2} >0$ and hence the best response set equals $\{0\}$. 
Further arguing as in Lemma \ref{lem_cont_of_attract}, one can show that best response set equals $\{0\}$ even at mutational strategy $\pi_{\epsilon}(\hat{\beta},p)$ with $\hbeta = 0$ when $\epsilon$ is small enough.  \eop }

\begin{lemma}\label{lem_cont_of_attract}
Let $\rho > 1$. Assume $\tilde{q}(\up)\ne 1$  where $q(\up) = \min\{\tilde{q}(\up), 1\}$. Consider a policy $\pi(\hat{\beta})$ where $\pi \in \Pi$ and $\hat{\up}$ is the attractor of the  corresponding   ODE   \eqref{eqn_ODE}. Let $\hat{\up}_\epsilon$ be attractor corresponding to $\epsilon$-mutant of this policy,  $\pi_\epsilon(\hat{\beta}, p)$ for some $p\in[0,1]$. Then, 
i) there exists an ${\bar \epsilon} (p) > 0$ such that the
attractor is unique  and  is a continuous function of $\epsilon $ 
for all $\epsilon \le {\bar \epsilon}$ with $\hat{\up}_0 = {\hat \up}$. 

\noindent ii) Further  $\bar{\epsilon}$ could be chosen such that the sign of $h(\hat{\up}_\epsilon)$ remains the same as that of $h(\hat{\up})$ for all $\epsilon \le {\bar \epsilon}$, when the latter is not zero. \eop
\end{lemma}

\noindent {\bf Proof:} We begin with an interior attractor. Such an attractor is a zero of a function like the  following (e.g.,  for  VFC1 it equals, see \eqref{eqn_ODE}): 

\vspace{-6mm}
\begin{align}\label{eqn_ODE_lemma}
\phi\lambda - r - b   ,\ \ \min \left \{1,  \psi \theta \hat{\beta} \right \} \phi \nu - b\psi , \mbox{ and } \frac{b-d }{\varrho} - \eta.  
\end{align}

\vspace{-2mm}
Under mutation policy, $\pi_\epsilon(\hat{\beta}, p)$, the function modifies to the following:

\vspace{-6mm}
\begin{align}\label{eqn_perturbed_ODE_lemma}
\phi\lambda - r - b   ,\ \ \Big( (1-\epsilon ) \mbox{$\min\{1,\hat{\beta} \psi \theta$\}} + \epsilon p\Big)  \phi \nu - b\psi , \mbox{ and } \frac{b-d }{\varrho} - \eta.  
\end{align}

\vspace{-2mm}
By directly computing the zero  of this function, it is clear that we again have unique zero  and these are continuous\footnote{When $\tilde{q}(\up) > 1$,  the zeros are
$\nicefrac{(\epsilon p + (1-\epsilon) )}{ (\mu \rho)}$, otherwise they are the zeros of a quadratic equation with varying parameters,  we have real zeros in this regime. 
} in $\epsilon$ (in some ${\bar \epsilon}$-neighbourhood) and that they coincide with $\hat{\up}$ at $\epsilon =0$. Further using Lyapunov function as defined in the corresponding proofs (with obvious modifications) one can show that these zeros are also attractors  in the neighborhood. 
\TR{The remaining part of the  proof is completed in \cite{TR}. }
{

Now consider the case when $\hat{\up} = (1-\nicefrac{1}{\rho}, 0)$ (the top row in all the attractor tables).  The first equation is not changed in \eqref{eqn_perturbed_ODE_lemma} from that in \eqref{eqn_ODE_lemma} whose zero provides $\hat{\theta}$ and the derivative near $\psi = 0 $ is negative\footnote{This is ensured by an appropriate $\epsilon$-neighborhood.} even for the second function defined in \eqref{eqn_perturbed_ODE_lemma} and the third function will have continuous zeros as before. Once again the Lyapunov arguments go through as in corresponding proofs. Thus in this case in fact
$\hat{\up}_\epsilon = {\hat \up}$  in the neighbourhood. 

The proof can be completed in exactly similar lines for the other attractors.
}
The last result follows by continuity of $h$ function \eqref{eqn_utility}. 
\eop



%
%

%
%
%

\hide{
\subsection*{To plot ODE solution versus simulated}
Start at $\bar K$ epochs and consider $M$ samples each at range $m$. 
Store simulated results at 
$$
\{\up^{sim}_{{\bar K} + mk} \}_{m\ge 1}
$$
with corresponding ODE time clocks 
$$
t_{{\bar K} + mk} 
:=  \sum_{i= {\bar K}}^{{\bar K} + mk} \epsilon_i - \epsilon_{\bar K},  \   \epsilon_i = 1/(i+1).
$$
Plot $\up^{sim}_{{\bar K} +mk}$ versus $t_{{\bar K} + mk} $.
Start ODE solution with initial condition $\up^{sim}_{{\bar K} }$.}

\end{document}